\documentclass{elsarticle}

\makeatletter
\def\ps@pprintTitle{
 \let\@oddhead\@empty
 \let\@evenhead\@empty
 \def\@oddfoot{\centerline{\thepage}}%
 \let\@evenfoot\@oddfoot}
\makeatother

\usepackage{amsmath,amsthm,amsfonts,amssymb,mathrsfs}
\usepackage{graphicx}
\usepackage[latin1]{inputenc}

\theoremstyle{definition}
\newtheorem{rem}{Remark}

\newcommand{\bR}{{\bf R}}
\newcommand{\bS}{{\bf S}}
\newcommand{\bT}{{\bf T}}
\newcommand{\slp}{{S_k}}
\newcommand{\Sw}{{S_k^{\rm \omega}}}
\newcommand{\Kw}{{K_k^{\rm \omega}}}
\newcommand{\Tw}{{T_k^{\rm \omega}}}

\newcommand{\dprime}{{T_k}}

\newcommand{\be}{\begin{equation}}
\newcommand{\ee}{\end{equation}}

\newcommand{\ba}{\begin{aligned}}
\newcommand{\ea}{\end{aligned}}

\def\lund{Centre for Mathematical Sciences, Lund University\\
  Box 118, 221 00 Lund, Sweden}

\def\ccm{Center for Computational Mathematics, Flatiron Institute, Simons Foundation\\
  New York, NY 10010, USA}

\begin{document}

\begin{frontmatter}

\title{The Helmholtz Dirichlet and Neumann problems on piecewise smooth open curves} 

\author[lund]{Johan Helsing}
\address[lund]{\lund}
\ead{johan.helsing@math.lth.se}

\author[ccm]{Shidong Jiang}
\address[ccm]{\ccm}
\ead{sjiang@flatironinstitute.org}

\begin{abstract}
  A numerical scheme is presented for solving the Helmholtz equation with Dirichlet or
  Neumann boundary conditions on piecewise smooth open curves, where the curves may have
  corners and multiple junctions. Existing integral equation methods for smooth open
  curves rely on analyzing the exact singularities of the density at endpoints for
  associated integral operators, explicitly extracting these singularities from the
  densities in the formulation, and using global quadrature to discretize the boundary
  integral equation. Extending these methods to handle curves with corners and multiple
  junctions is challenging because the singularity analysis becomes much more complex,
  and constructing high-order quadrature for discretizing layer potentials with singular
  and hypersingular kernels and singular densities is nontrivial.
  The proposed scheme is built upon the following two observations. First, the
  single-layer potential operator and the normal derivative of the double-layer
  potential operator serve as effective preconditioners for each other locally.
  Second, the recursively compressed inverse preconditioning (RCIP) method
  can be extended to
  address "implicit" second-kind integral equations. The scheme is high-order, adaptive,
  and capable of handling corners and multiple junctions without prior knowledge of
  the density singularity. It is also compatible with fast algorithms, such as the
  fast multipole method.
  The performance of the scheme is illustrated with several numerical examples.
\end{abstract}


\end{frontmatter}

\section{Introduction}
\label{sec:intro}
In this paper, we consider the time-harmonic wave scattering problem
from a collection of open curves in two dimensions. The governing
partial differential equation (PDE) for the scattered field $u$ is the
Helmholtz equation:
\be\label{eq:helmholtz}
\Delta u(r)+k^2 u(r)=0\,,
\quad r=(x,y)\in\mathbb{R}^2\setminus\Gamma\,,
\ee
with the Sommerfeld radiation condition
\be\label{eq:radiationbc}
\lim_{|r|\to\infty}
\sqrt{|r|}\left(\frac{\partial}{\partial|r|}-{\rm
  i}k\right)u(r)=0
\ee
imposed at infinity. Here, $k$ is the wavenumber, and
the boundary $\Gamma$ consists of a set
of piecewise smooth open curves. By ``piecewise smooth'', we mean that these curves
may have corners and branch points, such as triple junctions, quadruple
junctions, et cetera. On the boundary, we impose either the Dirichlet
condition
\be\label{eq:dirichlet}
u(r)=g(r)\,,\quad r\in\Gamma\,,
\ee
or the Neumann condition
\be\label{eq:neumann}
\frac{\partial u}{\partial\nu}(r)=g(r)\,,\quad r\in\Gamma\,,
\ee
where $\nu$ is the unit normal vector at $r$. These problems are frequently encountered
in acoustic and electromagnetic scattering. In acoustics, the Dirichlet condition
corresponds to a sound-soft object, and the Neumann condition corresponds to a
sound-hard object. In electromagnetics, they correspond to planar
transverse electric (TE) and transverse magnetic (TM) waves,
respectively~\cite{BrunLint12,helsing2024jmmct}.

Boundary integral methods are standard tools for solving the
Helmholtz Dirichlet and Neumann problems due to several advantages. First, the radiation
condition is satisfied automatically by the layer potential representation of the field,
removing the need of truncating the computational domain and constructing the so-called
artificial boundary conditions.
Second, the densities of layer potentials are supported on the boundary, leading to an
optimal number of unknowns during the solve phase. Third, the far-field pattern
can be evaluated accurately and efficiently via the integral representation and fast
algorithms, such as the fast multipole method (FMM)~\cite{wideband3d,wideband2d,greengard1998icse,rokhlin1993acha,rokhlin1998acha}. Among all possible integral formulations,
second-kind integral equations (SKIEs)
have become popular in the last two or three decades due
to their well-conditioning, which leads to stable and improved accuracy and
a low number of
iterations for iterative solvers such as GMRES~\cite{gmres}.

There are two main difficulties associated with piecewise smooth open curves. First,
there is a fundamental change in topology between closed curves and open curves. A
closed curve divides the entire space into two domains -- the interior and the exterior.
When the curve is open, there is no distinction between the interior and the exterior.
Thus, the so-called lightning solver~\cite{ginn2023arxiv}, a special case of method
of fundamental solutions~\cite{barnett2008jcp}, cannot be applied to solve these problems.
This is because the poles in those methods
must be placed outside the domain where the solution is sought, and the domain for the
open-curve problem is the entire space.
For closed curves,
the construction of SKIEs relies on the so-called jump relations of layer potentials
in potential theory.
For open curves, however, the field must be continuous across the boundary
for the Dirichlet problem.
And its normal derivative must be continuous across the boundary for
the Neumann problem. Thus, classical integral representations, which induce jumps across
the boundary, cannot be used for open curves. Furthermore,
the endpoints of an open curve introduce strong singularity in the density of layer
potentials, often changing the underlying function space and complicating the design
of numerical procedures.

Second, the singularity analysis of densities at corners and branch
points is rather difficult. Indeed, for the boundary integral
equation (BIE) formulation of the Helmholtz equation, the singularity
analysis at corners has been done only
recently~\cite{hoskins2019sisc,serkh2019acha,serkh2016jcp,serkh2016pnas}.
Since the analysis depends on
the governing equation, the boundary condition, and the type of point singularity,
it is carried out for each particular case. See, for example,
\cite{hoskins2020jcp} for the
Laplace equation with Neumann boundary conditions, \cite{rachh2020cpam} for the Stokes
equation with corners, and \cite{hoskins2020paa} for the Laplace equation with triple
junctions. The subsequent quadrature design, following the singularity analysis,
requires considerable effort. This is because the density often contains infinitely many
types of singularities that depend on the opening angle, say, in the case of corners.
Even with the use of generalized Gaussian
quadrature~\cite{ggq1,ggq2,ggq3}, a special quadrature must be constructed for each
specific case to ensure that the number of quadrature nodes is not excessively large.
The construction of these special quadratures is usually not
performed on the fly, as the precomputation of tables is costly
in terms of computation and storage.

In this paper, we propose a numerical scheme for solving the Helmholtz Dirichlet
and Neumann problems on open curves. The scheme is based on the following two
observations. First, 
the single layer potential operator $S_k$ (defined in \eqref{eq:Soper})
and the normal derivative of the double layer
potential operator $T_k$ (defined in \eqref{eq:Toper})
(or more generally, hypersingular operators including the quadruple
layer potential opeartor~\cite{kolm2003acha})
are good preconditioners of each other. For closed curves, these are
two well-known Calder{\' o}n identities from the Calder{\' o}n
projector~\cite{calderon1963,hsiao2008}. In \cite{costabel2003,stephan1987ieot},
it is shown that the densities
of the single-layer and hypersingular integral operators have inverse square root and
square root singularity in the vicinity of the endpoints, respectively. 
In \cite{BrunLint12}, the authors generalize these
Calder{\' o}n identities to the case of open curves by pulling out the aforementioned
singularities explicitly from the densities. However, it is nontrivial to extend
the work in \cite{BrunLint12} to treat the case of piecewise smooth open
curves, where the corners and branch points lead to additional difficulties.
In \cite{hang2009acha}, a {\it local} analysis is carried out to generalize
the Poincar{\' e}--Bertrand formula on a hypersurface. As the analysis is mostly local,
a similar argument can be used to show that $S_k$ and $T_k$ are good preconditioners
to each other locally. 
Second, the recursively compressed inverse preconditioning (RCIP)
method, first introduced in~\cite{helsing3}, has proven to be highly
successful in handling all types of point singularities for BIEs in two dimensions.
See~\cite{Tutorial} and
references therein. The RCIP method
provides a unified framework for addressing these challenging problems accurately and
efficiently, without requiring any analysis of the exact singularities in the densities
of layer potentials.

Our scheme combines these two observations in a straightforward manner. We apply
$T_k$ to precondition $S_k$ for the Dirichlet problem, and $S_k$ to precondition $T_k$
for the Neumann problem. Unlike \cite{BrunLint12}, we do not explicitly
extract the singularities from the densities. Although the authors of
\cite{BrunLint12} successfully extracted the singularities of the
densities for smooth open curves, doing so would be rather cumbersome numerically in
the presence of corners and branch points. Instead, we apply the RCIP method to
uniformly treat the endpoint and other point singularities.
We then divide the open curves into panels and apply a high-order kernel-split
quadrature~\cite{HelsKarl18,HelsOjal08} to discretize the layer potentials. It turns
out that the RCIP method is capable of figuring out the correct function spaces in a
remarkable way. Furthermore, even though the formulation is of the second-kind
{\it implicitly}, we have not observed any GMRES stagnation~\cite{AgaOneRac23} due to
the lack of numerical second-kindness in our numerical experiments.
The resulting numerical scheme is accurate, efficient, and fully
adaptive. It can also be easily coupled with fast algorithms, such as the FMM, to solve
large-scale problems. On the other hand, it seems rather difficult to extend the global
spectral discretization scheme presented in \cite{BrunLint12} to handle
corners and branch points, which may also have a compatibility issue with fast
algorithms~\cite{HaBaMaYo14}.

The rest of the paper is organized as follows. In Section~\ref{sec:preliminaries}, we
summarize the notations, analytical and numerical tools,
and the results in \cite{BrunLint12} for comparison purpose.
In Section~\ref{sec:mainresults}, we present our 
numerical scheme. Numerical examples are shown in
Section~\ref{sec:examples}.

\section{Preliminaries}\label{sec:preliminaries}

\subsection{Green's function and layer potentials}\label{sec:potentialtheory}

The free-space Green's function, also called fundamental solution, of 
the Helmholtz equation in $\mathbb{R}^2$ is
\begin{equation}
\Phi_k(r,r')=\frac{\rm i}{4}H_0^{(1)}(k|r-r'|)\,.
\label{eq:Phi}
\end{equation}
Here $k$ is a wavenumber, $r=(x,y)$, and $H_0^{(1)}$ is the zeroth
order Hankel function of the first kind~\cite{nisthandbook}.

Let $\Gamma$, as above, be a set of piecewise smooth open curves,
referred to collectively as the {\it boundary}. The unit normal at
$r\in\Gamma$ is $\nu(r)$. We use standard definitions of the Helmholtz
layer potential operators $S_k$, $K_k$, and
$T_k$~\cite[pp.~41--42]{ColtonKress98}, defined by their action on a
layer density $\rho(r)$ on $\Gamma$ as
\begin{align}
S_k\rho(r)&=2\int_{\Gamma}\Phi_k(r,r')\rho(r')\,{\rm d}\ell'\,, 
\quad r\in\Gamma\,,
\label{eq:Soper}\\
K_k\rho(r)&=
  2\int_{\Gamma}\frac{\partial\Phi_k}{\partial\nu'}(r,r')\rho(r')
            \,{\rm d}\ell'\,,\quad r\in\Gamma\,,
\label{eq:Koper}\\
T_k\rho(r)&=
  2\int_{\Gamma}\frac{\partial^2\Phi_k}
  {\partial\nu\partial\nu'}(r,r')\rho(r')\,{\rm d}\ell'\,,
\quad r\in\Gamma\,,
\label{eq:Toper}
\end{align}
where ${\rm d}\ell$ is an element of arc length,
$\partial/\partial\nu=\nu(r)\cdot\nabla$, and
$\partial/\partial\nu'=\nu(r')\cdot\nabla'$. For ease of presentation
we also use $S_k$ and $K_k$ to define single- and double-layer
potentials. That is, we use
\begin{align}
S_k\rho(r)&=2\int_{\Gamma}\Phi_k(r,r')\rho(r')\,{\rm d}\ell'\,, 
\quad r\in\mathbb{R}^2\setminus\Gamma\,,\\
K_k\rho(r)&=
  2\int_{\Gamma}\frac{\partial\Phi_k}{\partial\nu'}(r,r')\rho(r')
  \,{\rm d}\ell'\,,\quad r\in\mathbb{R}^2\setminus\Gamma\,.
\end{align}

\subsection{Calder{\' o}n identities on closed curves}
From~\cite[p.~43]{ColtonKress98} we have on closed curves
\begin{align}
  S_k(-T_k)&=I-K_kK_k\,,\quad r\in\Gamma\,,\label{eq:calderon1}\\
  T_k(-S_k)&=I-K_k^{\rm A}K_k^{\rm A}\,,\quad r\in\Gamma\,,\label{eq:calderon2}
\end{align}
where $I$ is the identity operator, and $K_k^{\rm A}$ is defined by
\be
K_k^{\rm A}\rho(r)=
2\int_{\Gamma}\frac{\partial\Phi_k}{\partial\nu}(r,r')\rho(r')
\,{\rm d}\ell'\,,\quad r\in\Gamma\,.
\label{eq:KAoper}
\ee

\subsection{Summary of the work in \cite{BrunLint12}}
Here, we provide a short summary of the work by Bruno and Lintner in \cite{BrunLint12}.
The main purpose is to make the comparison between our numerical scheme and the scheme
in \cite{BrunLint12} easy and transparent. Suppose now that $\Gamma$ is a smooth open
curve with the parameterization $r(t)=(x(t),y(t))$ for $t\in [-1,1]$. Define the weight
function $\omega$ by
\be
\omega(r(t))=\sqrt{1-t^2}.
\ee
Using the singularity analysis from \cite{costabel2003,stephan1987ieot}, Bruno and
Lintner propose to pull out the singularity from the density explicitly from $S_k$ and
$T_k$. To be more precise, they consider the weighted operators
\be\label{eq:weightedoperators}
\Sw[\rho](r) = \slp\left[\frac{\rho}{\omega}\right](r), \quad
\Tw[\rho](r) = \dprime[\omega\rho](r).
\ee
For the Dirichlet problem, the field is represented as
\be\label{eq:urepDB}
u(r) = \Sw[\rho](r), \quad  r\in\mathbb{R}^2\setminus\Gamma\,,
\ee
and the BIE is
\be\label{eq:DiriB}
\Tw\Sw\rho = \Tw g.
\ee
For the Neumann problem, the field is represented as
\be\label{eq:urepNB}
u(r) = \Kw\Sw[\rho](r), \quad  r\in\mathbb{R}^2\setminus\Gamma\,,
\ee
and the BIE is
\be\label{eq:NeumB}
\Tw\Sw\rho = g.
\ee
The BIEs (\ref{eq:DiriB},\ref{eq:NeumB})
 are further analyzed
with the cosine transform $t=\cos\theta$ 
and an even extension mapping the interval $[-1,1]$ to $[0,2\pi]$, so that
standard Calder{\' o}n identity \eqref{eq:calderon2} can be used. On the numerical
side, these BIEs are discretized
via a global spectral method using the fact that 
the convolution of the logarithmic kernel with trigonometric basis functions
can be evaluated analytically~\cite{kress2014}, which is equivalent to the fact
that the convolution integral of $\int_{-1}^1\ln|x-t| T_n(t)/\sqrt{1-t^2}dt$
can be evaluated analytically (see, for example, \cite[eqs. (27)-(28)]{JianRokh04}).
The authors further observe that the system matrix can be constructed
in $O(N^2\ln N)$ operations with the FFT, where $N$ is the total number of discretization
points on $\Gamma$.

The analysis in \cite{BrunLint12} and the subsequent paper
\cite{lintner2015prsea} is fairly complete, although there are some
drawbacks with the associated numerical scheme. First, it
is not clear how one could reduce the computational cost from $O(N^2\ln N)$ to
$O(N\ln N)$. A global kernel-split quadrature is applied to
discretize the associated integral operators. According to \cite{HaBaMaYo14}, the
quadrature scheme seems incompatible with the use of fast algorithms such as the FMM.
Second, being a global discretization scheme, it is inefficient to handle fine local
structure in the geometry and data. Third, there is no discussion on the evaluation
of the near field, though it is well-known that close evaluation is actually most
difficult in the design of high-order
quadratures for the evaluation of singular and nearly singular integrals.
Finally and most importantly, it is rather difficult to extend
the analysis and numerical scheme in \cite{BrunLint12} to treat open curves with
corners and branch points, as discussed in Section~\ref{sec:intro}.

\subsection{Overview of RCIP}
\label{sec:rcipoverview}

The RCIP method accelerates and stabilizes Nystr{\" o}m solvers in certain
``singular'' situations. More precisely, RCIP applies to the Nystr{\" o}m
discretization of Fredholm second kind BIEs on closed or open $\Gamma$
that contain some type of singular points. These points could be
endpoints, corners, or branch points of $\Gamma$. This section, based
on~\cite[Sections~3 and~16]{Tutorial}, reviews RCIP in a standard
setting. The Nystr{\" o}m discretization is
assumed to rely on composite Gauss--Legendre quadrature with $16$
nodes per quadrature panel. 

We start with a Fredholm second kind BIE on $\Gamma$ in standard form
\begin{equation}
  \left(I+G\right)\rho(r)=g(r)\,,\quad r\in \Gamma\,.
\label{eq:1}
\end{equation}
Here, $G$ is an integral operator that is often
assumed compact away from singular boundary points on $\Gamma$, $g(r)$
is a piecewise smooth right-hand side, and $\rho(r)$ is the unknown
layer density to be solved for.

For simplicity, we only discuss how to deal with one singular point on
$\Gamma$ and denote this point $\gamma$. Multiple singular points can
be treated independently and in parallel.

Let $G$ be split into two parts
\begin{equation}
G=G^\star+G^\circ\,,
\label{eq:splitG}
\end{equation}
where $G^\star$ describes the kernel interaction close to $\gamma$ and
$G^\circ$ is the (compact) remainder. Let the part of $\Gamma$ on
which $G^\star$ is nonzero be denoted $\Gamma^\star$.

Now introduce the {\it transformed density}
\begin{equation}
  \tilde{\rho}(r)=\left(I+G^\star\right)\rho(r)\,.
\label{eq:tilde}
\end{equation}
Then use~(\ref{eq:splitG}) and~(\ref{eq:tilde}) to
rewrite~(\ref{eq:1}) as
\begin{equation}
  \left(I+G^\circ(I+G^\star)^{-1}\right)\tilde{\rho}(r)
  =g(r)\,,\quad r\in \Gamma\,.
\label{eq:2}
\end{equation}
Although~(\ref{eq:2}) looks similar to~(\ref{eq:1}), there are
advantages to using~(\ref{eq:2}) rather than~(\ref{eq:1}) from a
numerical point of view, For example, the spectral properties of
$G^\circ(I+G^\star)^{-1}$ are better than those of $G$.

The Nystr{\" o}m/RCIP scheme discretizes~(\ref{eq:2}) chiefly on a grid on
a {\it coarse mesh} on $\Gamma$ that is sufficient to resolve
$G^\circ$, $\tilde{\rho}(r)$, and $g(r)$. Only the inverse
$(I+G^\star)^{-1}$ needs a grid on a locally refined {\it fine mesh}
on $\Gamma^\star$. The fine mesh is obtained from the coarse mesh by
$n_{\rm sub}$ times subdividing the coarse mesh in direction toward
$\gamma$. The discretization of~(\ref{eq:2}) assumes the form
\begin{equation}
  \left({\bf I}_{\rm coa}+{\bf G}_{\rm coa}^\circ{\bf R}\right)
  \tilde{\boldsymbol{\rho}}_{\rm coa}={\bf g}_{\rm coa}\,,
\label{eq:3}
\end{equation}
where ${\bf R}$ is a sparse block matrix called the {\it compressed
  inverse} and where the discrete unknown
$\tilde{\boldsymbol{\rho}}_{\rm coa}$ lives on the coarse grid only.
The boundary part $\Gamma^\star$ is assumed to contain the two, four,
or six coarse panels closest to $\gamma$ depending on if $\gamma$ is
an endpoint, a corner, or a branch point.

The power of RCIP lies in the construction of ${\bf R}$. In theory,
${\bf R}$ corresponds to a discretization of $(I+G^\star)^{-1}$ on the
fine grid, followed by a lossless compression back to the coarse grid.
In practice, ${\bf R}$ is constructed via a forward recursion of
length $n_{\rm sub}$ on a hierarchy of intermediate grids on
$\Gamma^\star$, and where refinement and compression occur in
tandem. The computational cost grows, at most, linearly with
$n_{\rm sub}$.

Once the compressed equation~(\ref{eq:3}) is solved for
$\tilde{\boldsymbol{\rho}}_{\rm coa}$ and its {\it weight-corrected}
counterpart
\begin{equation}
\hat{\boldsymbol{\rho}}_{\rm coa}
\equiv{\bf R}\tilde{\boldsymbol{\rho}}_{\rm coa}
\label{eq:wcorr0}
\end{equation}
is produced, several useful functionals of $\rho(r)$ can be computed
with ease. Should one so wish, the solution $\rho(r)$ on the fine grid
can be reconstructed from $\tilde{\boldsymbol{\rho}}_{\rm coa}$ via a
backward recursion on the same hierarchy of grids that is used in the
construction of ${\bf R}$.

\begin{rem}
  Note that the number of steps in the recursion for ${\bf R}$ is the
  same as the number of subdivisions needed to create the fine mesh
  from the coarse mesh and that this number is denoted $n_{\rm sub}$.
  In the numerical experiments of Section~\ref{sec:examples} we shall
  illustrate the convergence of a field quantity as a function of
  $n_{\rm sub}$.
\end{rem}
  
\begin{rem}
  The assumption, mentioned above, that $G^\circ$ of~(\ref{eq:splitG})
  should be compact for RCIP to apply, is actually too
  restrictive. See~\cite[Eq.~(24)]{Tutorial} for a better, although
  more technical, assumption about $G^\circ$. Allowing $G^\circ$
  of~(\ref{eq:splitG}) to be non-compact is important in the present
  work. In~(\ref{eq:AB}) we view $B$ as a reasonable right
  preconditioner to $A$ and merely consider $(A(-B)-I)^\circ$ as an
  operator whose spectrum is clustered.
\end{rem}

\subsection{On-surface and near-field evaluation}
\label{sec:closeevaluation}

The discretization of our BIEs on $\Gamma$ involves the discretization
of integrals with various types of singular integrands. Furthermore,
for $r$ close to $\Gamma$, the discretization of the field
representations involves various types of nearly singular
integrands. These singular and nearly singular integrals can,
typically, not be accurately evaluated at a reasonable cost using
composite standard interpolatory quadrature, but require special
discretization techniques. For this, we use a panel-based product
integration scheme which, according to the classification
of~\cite{HaBaMaYo14}, is ``explicit split''.

Our panel-based explicit-split quadrature scheme is described in
detail in~\cite[Section~4]{HelsKarl18}. For now, we only mention that
the scheme combines recursion with analytical methods such as the use
of the fundamental theorem of complex line integrals. Nearly singular
and singular integral operator kernels $G(r,r')$ are split according
to the general pattern
\begin{multline}
G(r,r')\,{\rm d}\ell'=G_0(r,r')\,{\rm d}\ell'
+\log|r-r'|G_{\rm L}(r,r')\,{\rm d}\ell'\\
+c_{\rm C}\Re{\rm e}\left\{\frac{G_{\rm C}(z,t)\,{\rm d}t}
                 {{\rm i}(t-z)}\right\}
+\Re{\rm e}\left\{\frac{G_{\rm H}(z,t)\,{\rm d}t}
                 {{\rm i}(t-z)^2}\right\}\,.
\label{eq:ex3}
\end{multline}
Here $c_{\rm C}$ is a, possibly complex, constant and $G_0(r,r')$,
$G_{\rm L}(r,r')$, $G_{\rm C}(z,t)$, and $G_{\rm H}(z,t)$ are smooth
functions. Our scheme requires explicit formulas for $G(r,r')$,
$G_{\rm L}(r,r')$, $G_{\rm C}(z,t)$, and $G_{\rm H}(z,t)$, while
$G_0(r,r')$ needs only to be known in the limit $r'\to
r$. See~\cite[Section~6.2]{Ande23} for additional remarks.

\section{The new numerical scheme}\label{sec:mainresults}
\subsection{Integral representations and boundary integral equations}
For the Dirichlet problem, we represent the field
$u(r)$ via
\begin{equation}
  u(r)=S_k(-T_k)\rho(r)\,,\quad r\in\mathbb{R}^2\,.
\label{eq:urepD}
\end{equation}
The BIE is
\begin{equation}
  S_k(-T_k)\rho(r)=g(r)\,,\quad r\in\Gamma\,.
\label{eq:Diri}
\end{equation}
For the Neumann problem, the field is represented via
\begin{equation}
  u(r)=K_k(-S_k)\rho(r)\,,\quad r\in\mathbb{R}^2\setminus\Gamma\,,
\label{eq:urepN}
\end{equation}
and the BIE is
\begin{equation}
  T_k(-S_k)\rho(r)=g(r)\,,\quad r\in\Gamma\,.
\label{eq:Neum}
\end{equation}
We would like to comment on some differences between the BIEs
(\ref{eq:Diri},\ref{eq:Neum})
and (\ref{eq:DiriB},\ref{eq:NeumB}) in \cite{BrunLint12}.
First, the hypersingular integral operator $T_k$ is used as the right preconditioner
in our representation. The modification is somewhat trivial, but
it removes the need to evaluating the action of $T_k$ on the boundary data $g$
in \cite{BrunLint12}. This is important when the open curves contain corners
and branch points, as otherwise one would need to apply proper weight functions in
$T_k$ to ensure that the modified right-hand side is piecewise smooth.
Second, the weight function $\omega$ is completely removed from the integral
representations (\ref{eq:urepD},\ref{eq:urepN}). At the first sight, this may seem like
a setback, as the strong inverse square root singularity changes the underlying
function space for the density, and the BIEs (\ref{eq:Diri},\ref{eq:Neum})
are no longer of the second kind.

However, as observed in \cite{hang2009acha}, $S_k$ and $T_k$ are good preconditioners
for each other pointwise, except at singular points. Thus, the BIEs
(\ref{eq:Diri},\ref{eq:Neum}) will be well-conditioned if the density singularities
in the vicinity of those singular points (that is, endpoints, corners, and branch points)
can be treated effectively. This is where the RCIP method comes into play.

\subsection{Incorporating the RCIP into the BIEs (\ref{eq:Diri},\ref{eq:Neum})}
\label{sec:compose2}

The BIEs~(\ref{eq:Diri}) and (\ref{eq:Neum}) are cast in the form
\begin{equation}
  A(-B)\rho_1(r)=g(r)\,,\quad r\in \Gamma\,,
\label{eq:AB}
\end{equation}
where the unknown density now is denoted $\rho_1(r)$ and where the
integral operators $A$ and $B$ are either $S_k$ and $T_k$ or $T_k$ and
$S_k$. Clearly, the composition of $A$ and $-B$ and the absence of the
identity operator
make~(\ref{eq:AB}) look quite different from the standard
form~(\ref{eq:1}).

We now apply the Nystr{\" o}m/RCIP scheme of Section~\ref{sec:rcipoverview} to
(\ref{eq:AB}), with the goal to find the analogue of the compressed
discrete equation~(\ref{eq:3}) for~(\ref{eq:AB}). First, introduce a
temporary density $\rho_2(r)$ via
\begin{equation}
\rho_2(r)=-B\rho_1(r)\,.
\end{equation}
This allows us to rewrite~(\ref{eq:AB}) as the $2\times 2$ block
system
\begin{equation}
\left(
\begin{bmatrix}
I & 0 \\
0 & I
\end{bmatrix} 
+
\begin{bmatrix}
 -I & A \\
  B & 0
\end{bmatrix} 
\right)
\begin{bmatrix}
 \rho_1(r) \\ 
 \rho_2(r)
\end{bmatrix}
=
\begin{bmatrix}
   g(r) \\ 0 
\end{bmatrix},
\label{eq:11}
\end{equation}
which corresponds to~(\ref{eq:1}) in Section~\ref{sec:rcipoverview}.

Note that~(\ref{eq:11}) is free from composed operators so that the standard
RCIP method can be readily applied. The analogue
of~(\ref{eq:2}) is
\begin{equation}
\left(
\begin{bmatrix}
I & 0 \\
0 & I
\end{bmatrix} 
+
\begin{bmatrix}
 -I^\circ & A^\circ \\
  B^\circ & 0
\end{bmatrix} 
\left(
\begin{bmatrix}
I & 0 \\
0 & I
\end{bmatrix} 
+
\begin{bmatrix}
 -I^\star & A^\star \\
  B^\star & 0
\end{bmatrix} 
\right)^{-1}\right)
\begin{bmatrix}
 \tilde{\rho}_1(r) \\ 
 \tilde{\rho}_2(r)
\end{bmatrix}
=
\begin{bmatrix}
   g(r) \\ 0 
\end{bmatrix}.
\label{eq:22}
\end{equation}
The analogue of~(\ref{eq:3}) is
\begin{equation}
\left(
\begin{bmatrix}
{\bf I}_{\rm coa} & {\bf 0}_{\rm coa} \\
{\bf 0}_{\rm coa} & {\bf I}_{\rm coa}
\end{bmatrix} 
+
\begin{bmatrix}
 -{\bf I}^\circ_{\rm coa} & {\bf A}^\circ_{\rm coa} \\
  {\bf B}^\circ_{\rm coa} & {\bf 0}_{\rm coa}
\end{bmatrix} 
\begin{bmatrix}
{\bf R}_1 & {\bf R}_3 \\
{\bf R}_2 & {\bf R}_4
\end{bmatrix} 
\right)
\begin{bmatrix}
 \tilde{\boldsymbol{\rho}}_{1{\rm coa}} \\ 
 \tilde{\boldsymbol{\rho}}_{2{\rm coa}}
\end{bmatrix}
=
\begin{bmatrix}
   {\bf g}_{\rm coa} \\ {\bf 0}_{\rm coa}
\end{bmatrix}.
\label{eq:33}
\end{equation}
Here we have partitioned the compressed inverse ${\bf R}$ into the
four sparse equi-sized blocks ${\bf R}_1$, ${\bf R}_2$ ${\bf R}_3$,
and ${\bf R}_4={\bf I}^\circ_{\rm coa}$.

As an extra twist, we write the compressed equation~(\ref{eq:33}) as a
linear system involving only one discrete density -- not two. To this
end, introduce the new density $\tilde{\boldsymbol{\rho}}_{\rm coa}$
via
\begin{equation}
\begin{split}
  \tilde{\boldsymbol{\rho}}_{1{\rm coa}}&=
  \tilde{\boldsymbol{\rho}}_{\rm coa}+
  {\bf R}_1^{-1}{\bf R}_3{\bf B}_{\rm coa}^\circ{\bf R}_1
  \tilde{\boldsymbol{\rho}}_{\rm coa}\,,\\
  \tilde{\boldsymbol{\rho}}_{2{\rm coa}}&=
  -{\bf B}_{\rm coa}^\circ{\bf R}_1
  \tilde{\boldsymbol{\rho}}_{\rm coa}\,.
\end{split}
\label{eq:single12}
\end{equation}
The change of variables~(\ref{eq:single12}) is chosen so that the
second block-row equation of~(\ref{eq:33}) is automatically
satisfied. The first block-row equation of~(\ref{eq:33}) becomes
\begin{multline}
  \left({\bf I}_{\rm coa}^\star
    -{\bf A}_{\rm coa}^\circ
\left({\bf R}_4-{\bf R}_2{\bf R}_1^{-1}{\bf R}_3\right)
{\bf B}_{\rm coa}^\circ{\bf R}_1+{\bf A}_{\rm coa}^\circ{\bf R}_2
\right.\\
\left.+{\bf R}_1^{-1}{\bf R}_3{\bf B}_{\rm coa}^\circ{\bf R}_1
\right)\tilde{\boldsymbol{\rho}}_{\rm coa}={\bf g}_{\rm coa}\,.
\label{eq:disc3}
\end{multline}

We observe, from~(\ref{eq:single12}), that away from $\Gamma^\star$
the discrete density $\tilde{\boldsymbol{\rho}}_{\rm coa}$ coincides
with $\tilde{\boldsymbol{\rho}}_{1{\rm coa}}$. Therefore one can
expect the system~(\ref{eq:disc3}) to share more properties with the
original equation~(\ref{eq:AB}) than the expanded system~(\ref{eq:33})
does. Numerical experiments indicate that~(\ref{eq:disc3}) is superior
to~(\ref{eq:33}) not only in terms of computational economy, but also
in terms of stability and of convergence of iterative solvers.

The analogue of~(\ref{eq:wcorr0}) is
\begin{equation}
\begin{split}
\hat{\boldsymbol{\rho}}_{1{\rm coa}}
  &={\bf R}_1\tilde{\boldsymbol{\rho}}_{\rm coa}\,,\\
\hat{\boldsymbol{\rho}}_{2{\rm coa}}
  &={\bf R}_2\tilde{\boldsymbol{\rho}}_{\rm coa}
    -\left({\bf R}_4-{\bf R}_2{\bf R}_1^{-1}{\bf R}_3\right)
    {\bf B}^\circ_{\rm coa}{\bf R}_1\tilde{\boldsymbol{\rho}}_{\rm coa}\,.
\end{split}
\end{equation}

\subsection{Some algorithmic details}
As mentioned in Section~\ref{sec:rcipoverview}, Nystr{\" o}m discretization is
used for integral equations. It relies solely on composite $16$-point
Gauss--Legendre quadrature and product integration based on polynomial
interpolation as described in Section~\ref{sec:closeevaluation}.
We use the low-threshold stagnation-avoiding GMRES of
\cite[Section~8]{HelsOjal08} to solve the resulting linear system,
with the matrix-vector product accelerated
by the FMM. The FMM computes the sum of the following form
\be
u_i=\sum_{j=1}^N K(r_i,r'_j)q_j, \quad j=1,\cdots,N, \quad i=1,\cdots,M\,,
\ee
where $K$ is a kernel function,
$M$ is the number of target points,
and $N$ is the number of source points. For highly oscillatory
Helmholtz kernels, the cost is $O((M+N)\log(M+N))$~\cite{wideband2d,greengard1998icse}.
When the set of target points is identical to the set of source
points, the self-interaction terms $j=i$ are excluded as most
kernels are singular at the origin.

We use the single-layer potential operator $S_k$ as an illustration
of computational complexity in the context of solving discretized
integral equations. To avoid cumbersome indices, 
we denote the matrix discretization of the operator 
$S_k^\circ$ simply as $\bS^\circ$.
In order to apply
$\bS^\circ$ fast to a given vector, we split it into two parts:
\be
\bS^\circ = \bS^\circ_{\rm smooth}+\bS^\circ_{\rm qc},
\ee
where $\bS^\circ_{\rm smooth}$ is obtained by discretizing $S_k^\circ$ via standard
Gauss-Legendre quadrature on each panel, and
$\bS^\circ_{\rm qc}$ is the quadrature correction part due to the logarithmic singularity
of the kernel. We further write $\bS^\circ_{\rm smooth}$ as
\be
\bS^\circ_{\rm smooth}=\bS_{\rm smooth}-\bS^\star_{\rm smooth}.
\ee
Since the weights are functions of the
source point $r'$ only, it is clear that the action of $\bS_{\rm smooth}$ can be taken
by the FMM. The quadrature correction part $\bS^\circ_{\rm qc}$ depends on both the target
and source points, but it is nonzero only for target points close to the panel for each
souce panel. Furthermore, the bad part $\bS^\star_{\rm smooth}$ is nonzero for panels
that are near the singular points. Thus, both $\bS^\circ_{\rm qc}$
and $\bS^\star_{\rm smooth}$ are sparse and can be constructed and
applied in $O(N)$ time and storage for most
geometries. Combining all this, we observe that $\bS^\circ$ can be applied in
$O(N\log N)$ time.
The preconditioner $\bR$ is the same as the identity
matrix, except for the panels close to the singular points. The cost of constructing
the nontrivial blocks of $\bR$ for each singular point is
$n_{\rm edge}^3 n_{\rm gl}^3 n_{\rm sub}$, where $n_{\rm edge}$ is the number of edges connected
to the singular point, $n_{\rm gl}$ is the number of Gauss-Legendre nodes on each panel,
and $n_{\rm sub}$ is the level of dyadic refinements along each edge in the
forward recursion. Note that $\bR$
can be reused for identical geometries, for example, corners with the same opening
angle, which is not uncommon in applications.

The evaluation of the field follows a similar approach. We decompose the operator
into a smooth component and a quadrature correction component. The smooth part
is computed using the FMM with $O((N+M)\log(N+M))$ cost,
while the quadrature correction is computed directly.
Since only $O(1)$ nearby targets require quadrature correction for each source panel,
the cost of applying the correction is $O(N+M)$ for typical applications.
It is important to note that a sorting algorithm is needed to identify the targets
requiring quadrature correction for each source panel. This can be done efficiently
using either an adaptive quadtree~\cite{biros2008sisc} or a $k$-d tree
(see, for example, \cite{saye2014camcs,yu2015acm}).

\section{Numerical examples}\label{sec:examples}
The numerical codes used are implemented in {\sc Matlab}, release
2022a, and executed on
a laptop equipped with eight Intel i9-10885H CPU Cores. 
We use the wideband Helmholtz FMM from the {\tt fmm2d}~\cite{fmm2d} library,
compiled with OpenMP enabled; the parallel efficiency ranges from 50\% to 70\%.

In most examples, the total field is given by
$u^{\rm tot}(r)=u^{\rm in}(r)+u^{\rm sc}(r)$,
where $u^{\rm in}$ represents an incident plane wave
\be\label{eq:incfield}
u^{\rm in}(r) = e^{i k (x\cos\theta+y\sin\theta)}
\ee
with the incident angle $\theta$, and
$u^{\rm sc}$ denotes the scattered field. We assume that
the total field satifies a zero boundary condition,
which leads to an appropriate boundary condition for the scattered field,
solved numerically.
We begin with a coarse mesh on $\Gamma$, which is
sufficient to resolve $G^\circ$, $\tilde{\rho}(r)$, and $g(r)$ to the
prescribed tolerance. For field and error images in Figures~\ref{fig:spiral_field},
\ref{fig:corner_field}, \ref{fig:branch_field}, and \ref{fig:maze_field},
GMRES tolerance is set to $10^{-13}$.
When assessing the accuracy of computed fields we
adopt a procedure where, for each numerical solution, we also compute
an overresolved reference solution, using approximately 50\% more points in the
discretization of the setup under study.
Both the fields and their errors are
computed at $3000\times 3000$ points on a Cartesian grid
in the computational domains shown. The absolute difference
between these two solutions at a field point is referred to as the {\it
  estimated absolute pointwise error}.

\subsection{Validation of the correctness of formulation and implementation}
\label{sec:validation}

As there are no simple analytic solutions to the Helmholtz Dirichlet and Neumann
problems on open curves, we use a Fortran code that implements the
$u(r) = S_kS^{-1}_{0,I}\rho(r)$ formulation for the Dirichlet problem
(a straightforward extension
of the formulation in \cite{JianRokh04} and also used in \cite{BrunLint12}
as a comparison) to cross validate the correctness of the numerical scheme
and its implementation in this paper. Here, $S_{0,I}^{-1}$ is the inverse
of the Laplace single layer potential operator on the straight line segment,
which can be computed analytically for Chebyshev polynomials.

In this first example, $\Gamma$ is the straight line segment
parameterized as
\begin{equation}
r(s)=(s,-0.2)\,, \quad s\in[-1,1]\,.
\label{eq:straight}
\end{equation}  
The wavenumber is $k=3$, the boundary data is
\begin{equation}
  g(r)=4x^3+2x^2-3x-1\,,
\label{eq:gstraight}
\end{equation}
and we seek $u(r)$ at a target point $r_{\rm targ}=(0.17,0.62)$.
The reference value
$$u_{\rm ref}=0.02788626934981090-0.75932847390327920{\rm i}$$
for $u(r_{\rm targ})$ is
obtained using $24$-digit variable-precision arithmetic and $24$-point
composite Gauss--Legendre quadrature in the discretization of the
formulation (\ref{eq:urepD},\ref{eq:Diri}).
The value produced by the Fortran code mentioned above agrees with
$u_{\rm ref}$ to about 14 digits, thereby validating the correctness
of both the formulation (\ref{eq:urepD},\ref{eq:Diri}) itself and
of our Nystr{\" o}m/RCIP discretization thereof.

\begin{figure}[!ht]
\centering
\includegraphics[height=45mm]{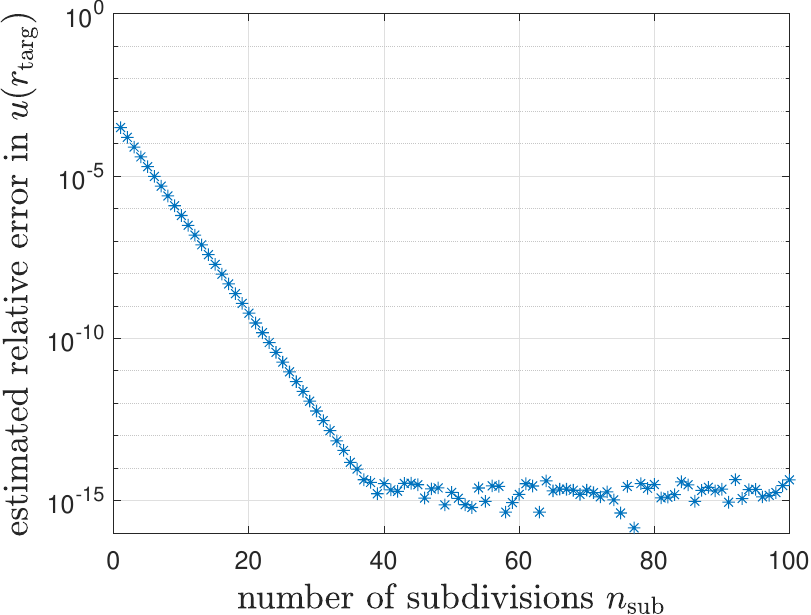}
\caption{\sf Convergence of $u(r_{\rm targ})$ with $n_{\rm sub}$ for the
  Helmholtz Dirichlet problem on the straight line segment~(\ref{eq:straight}).} 
\label{fig:line_nsub}
\end{figure}

Like all RCIP-based numerical schemes, the convergence of field
quantities as a function of $n_{\rm sub}$ is an important measure on
the stability and accuracy of our numerical scheme. 
Figure~\ref{fig:line_nsub} shows that the error decreases to about
$2\times 10^{-15}$ as $n_{\rm sub}$ increases from $1$ to $39$, then fluctuates
below $10^{-14}$ afterwards all the way to $n_{\rm sub}=100$.
Here the coarse mesh on $\Gamma$ has $6$ panels, that is, $96$ discretization points.
An estimated relative residual of machine epsilon ($\epsilon_{\rm mach}$)
in double-precision floating point arithmetic is used as the
GMRES stopping criterion (tolerance). The number of GMRES iterations needed for
full convergence with the formulation (\ref{eq:urepD},\ref{eq:Diri})
is $17$. Since $n_{\rm sub}$ is the level of dyadic refinement
towards the endpoints and $\epsilon_{\rm mach}=2^{-52}$,
Figure~\ref{fig:line_nsub} clearly shows that our numerical scheme is
as stable and accurate as the standard RCIP method applied to
an explicit SKIE, even though a hypersingular integral operator is used
in the formulation. Furthermore, there is no GMRES stagnation, even though
the formulation does not possess the so-called numerical second-kindness
in \cite{AgaOneRac23}.

\subsection{Comparison with the results in \cite{BrunLint12}}
In this section, we present conditioning and timing results on a
straight line segment and a spiral that are directly comparable with
\cite[Tables~2--5]{BrunLint12}. The total number of discretization
points on $\Gamma$ is $N$.

In \cite{BrunLint12} a global spectral
method is used to discretize SKIEs (\ref{eq:DiriB},\ref{eq:NeumB}),
the whole system matrix is built using FFT in $O(N^2\log N)$ time,
and the resulting linear system is then solved using GMRES. The prescribed precision
is $10^{-5}$. However, it is not clear whether the field is also computed to precision
$10^{-5}$, since there is no mentioning about the special quadrature
on the close evaluation of the field and the number of points used in generating field
plots \cite[Figure~3]{BrunLint12}. The timing is measured via {\sc Matlab}'s
command pair {\tt tic} and {\tt toc}.

In our numerical experiments, the GMRES tolerance is set to $10^{-6}$, the field
is evaluated at $300\times 300$ equispaced tensor grid on a rectangle. We vary
the wavenumber $k$ so that the ratio $L/\lambda = 50 \cdot 2^{i}$ for $i=0,1,\ldots,9$,
where $L$ is the arc length of $\Gamma$ and $\lambda=2\pi/k$ is the
wavelength. It is clear that the cases $L/\lambda=50, 200, 800$
presented in \cite[Tables~2--5]{BrunLint12} is a subset of what we study here.
We estimate the relative $l_2$ error by comparing with a reference
solution obtained by increasing $N$ and $n_{\rm sub}$ by about
$50\%$ each. That is,
\be
E_2 = \frac{\|u_{\rm num}-u_{\rm ref}\|_2}{\|u_{\rm ref}\|_2}.
\ee

\subsubsection{Straight line segment}

\begin{figure}[t]
\centering
\includegraphics[height=45mm]{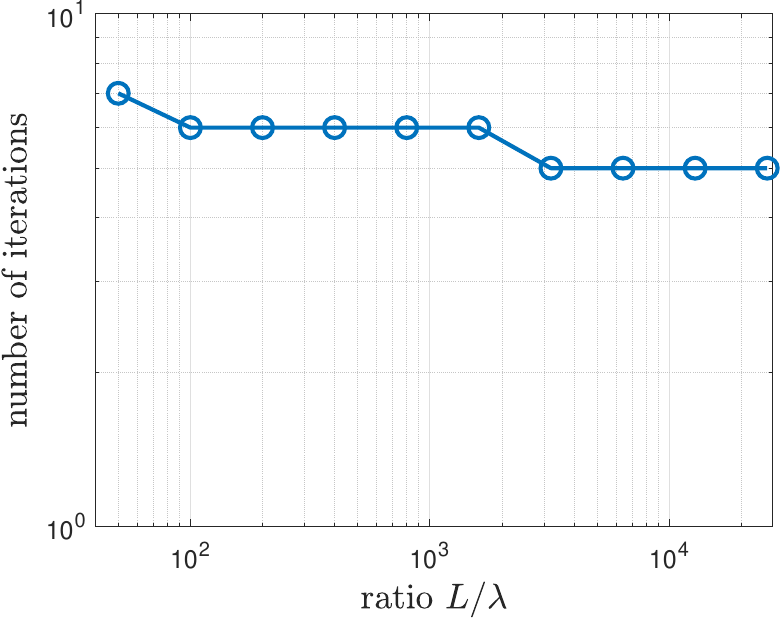}
\hspace*{2mm}
\includegraphics[height=45mm]{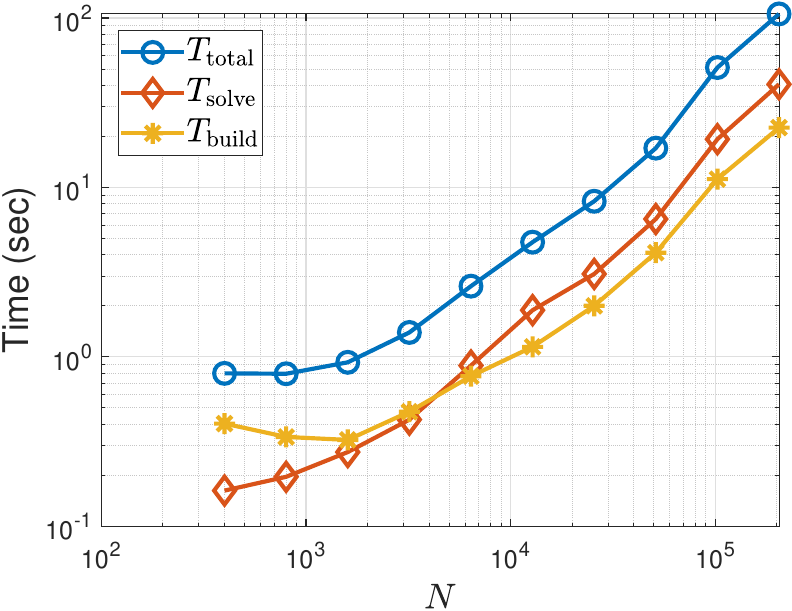}\\

\vspace{2mm}

\includegraphics[height=45mm]{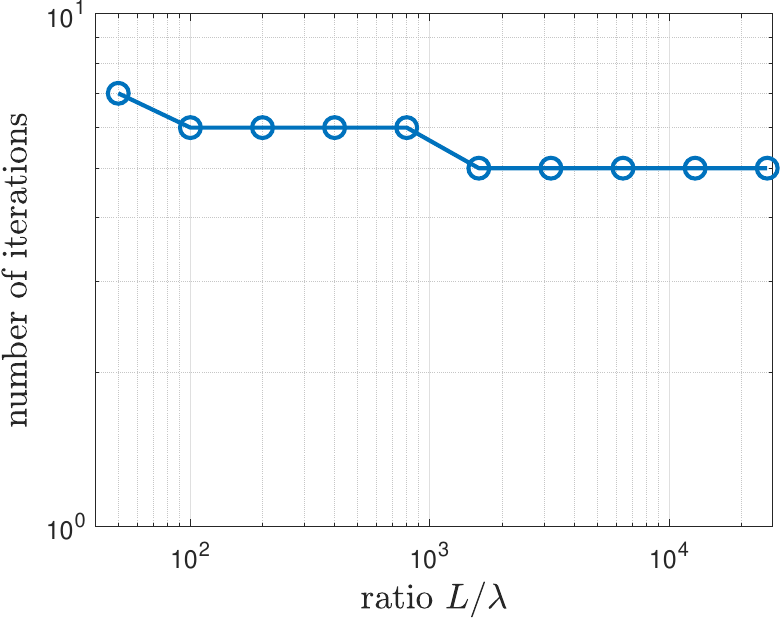}
\hspace*{2mm}
\includegraphics[height=45mm]{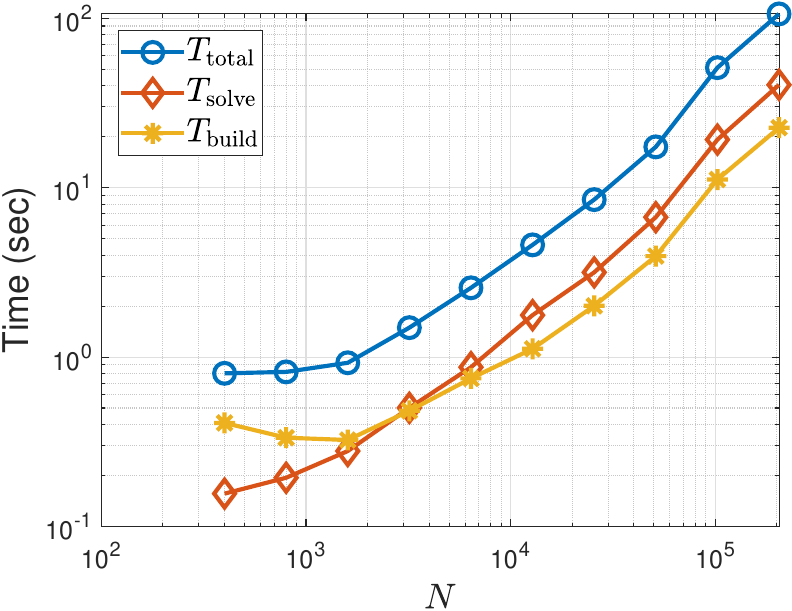}\\
\caption{\sf Number of GMRES iterations and computation times as a function of
wavenumber for the Helmholtz problem on a straight line segment.
Top: Dirichlet condition. Bottom: Neumann condition. Left: Number of GMRES iterations
as a function of $L/\lambda$. Right: timing results as a function of $N$ -- the total
number of discretization points on $\Gamma$. Here, $\Gamma$ is divided into
equi-sized panels in the parameter space and $16$ Gauss-Legendre
nodes are used to discretize each panel. $N=8L/\lambda$ for this example.
The relative $l_2$ error is measured on a
$300\times 300$ tensor grid over the square $[-1.3, 1.3]\times[-1.5, 1.1]$. The largest
relative $l_2$ errors are about $4\times 10^{-7}$ for the Dirichlet problem
and $6\times 10^{-7}$ for the Neumann problem.
}
\label{fig:line_iter}
\end{figure}

Figure~\ref{fig:line_iter} shows the conditioning and timings of our solver
for the Helmholtz Dirichlet and Neumann problems on the straight
line segment~\eqref{eq:straight}, compare \cite[Table~2]{BrunLint12} for
the Dirichlet problem and \cite[Table~4]{BrunLint12} for the Neumann problem.
Here, $T_{\rm total}$ is the total computation time; $T_{\rm build}$ is the
time on building the quadrature part of the system matrix $\bS^\circ_{\rm qc}$
and $\bT^\circ_{\rm qc}$, the smooth part of the system matrix $\bS^\star_{\rm smooth}$
and $\bT^\star_{\rm smooth}$ in the vicinity of the singular points, and the
preconditioner $\bR$; and $T_{\rm solve}$ is the time on GMRES.
We observe that
the Dirichlet problem is very similar to the Neumann problem in terms of both
conditioning and timings. Amazingly, the number of points for achieving six digits
of accuracy for the straight line segment is exactly the same as that
in \cite{BrunLint12} for achieving five digits of accuracy, even though we are using
a panelwise discretization, while a global spectral method is used in \cite{BrunLint12}.
In \cite{BrunLint12}, the number of GMRES iterations is eight for
the Dirichlet problem and nine for the Neumann problem, while one
only needs at most seven iterations, and often only five or six, for
both the Dirichlet and the Neumann problem with our scheme. This
difference is certainly very minor, but it does demonstrate that our
more general and purely numerical scheme performs as well as the
analysis-based scheme in \cite{BrunLint12} in an example where the
latter scheme is applicable.

For timings, both $T_{\rm build}$ and $T_{\rm solve}$ in our numerical
methods are less than one second for $L/\lambda=800$, while it takes more than
$50$ seconds for building the matrix and $14$ seconds for solving the linear system
in \cite{BrunLint12}.
The difference will be larger for smaller wavelengths, requiring
larger $N$, since our scheme scales like $O(N)$
for $T_{\rm build}$, and $O(N_{\rm iter}N\log N)$ for $T_{\rm solve}$, with $N_{\rm iter}$
the number of iterations in GMRES, while the scheme in \cite{BrunLint12} scales like
$O(N^2 \log N)$ for $T_{\rm build}$ and $O(N_{\rm iter} N^2)$ for $T_{\rm solve}$. 
We set $n_{\rm sub}=40$ for computing the numerical solution and $n_{\rm sub}=60$
for computing the reference solution. The time in computing $\bR$ is about $0.2$ seconds
for two endpoints.

\begin{figure}[t]
\centering
\includegraphics[height=45mm]{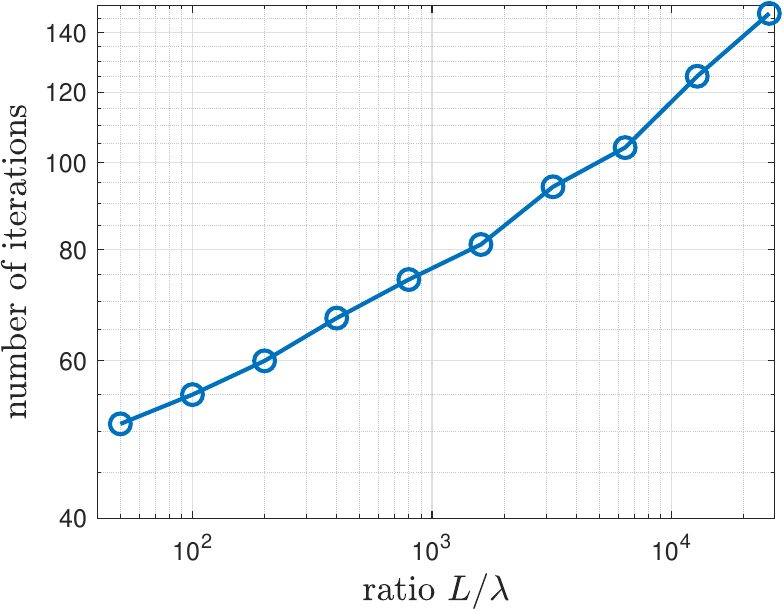}
\hspace*{2mm}
\includegraphics[height=45mm]{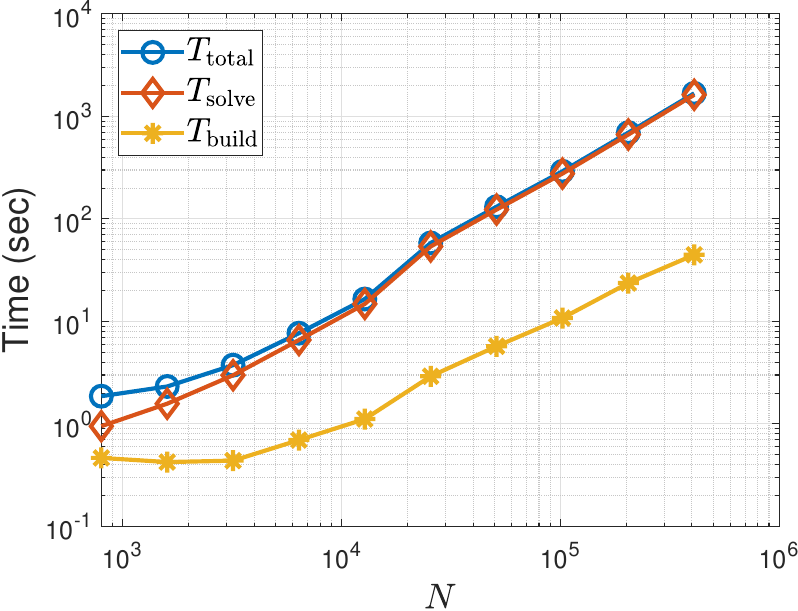}\\

\vspace{2mm}

\includegraphics[height=45mm]{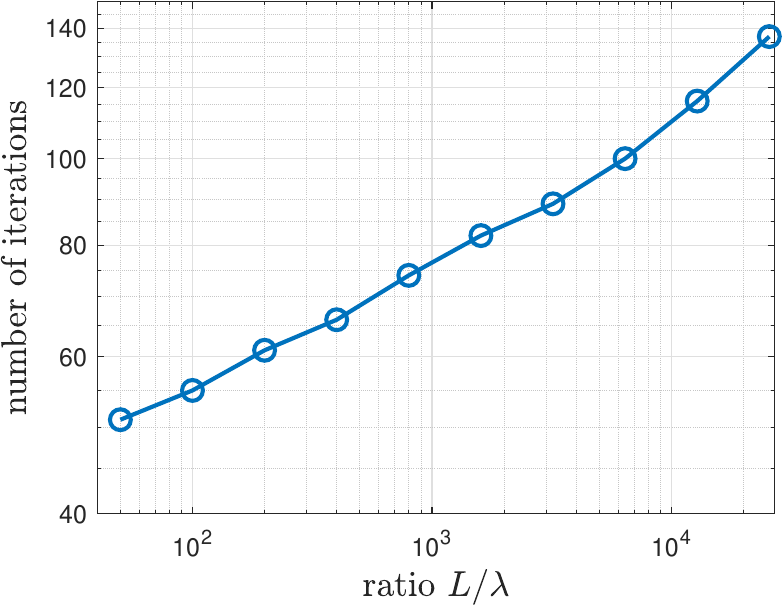}
\hspace*{2mm}
\includegraphics[height=45mm]{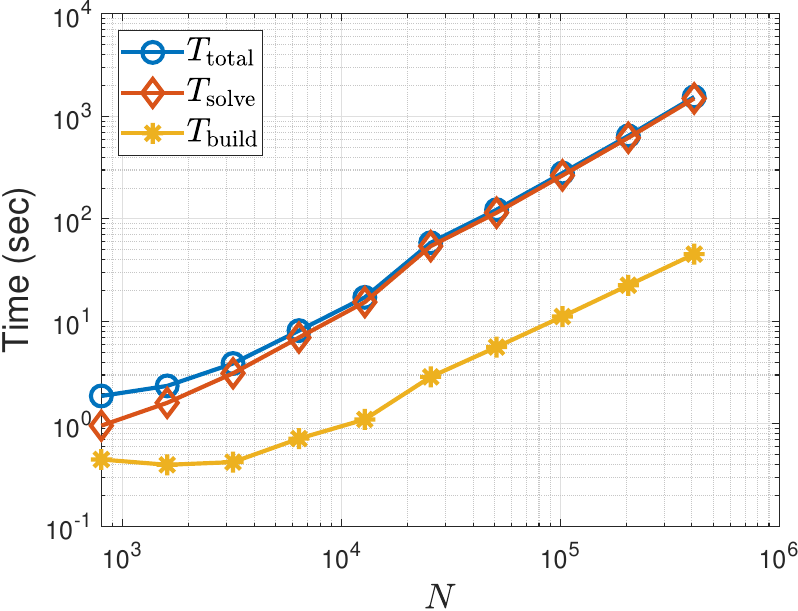}\\
\caption{\sf Same as Figure~\ref{fig:line_iter}, but for 
the spiral defined by \eqref{eq:spiral}.
$N=16L/\lambda$ for this example.
}
\label{fig:spiral_iter}
\end{figure}

\begin{figure}[!ht]
\centering
\includegraphics[height=45mm]{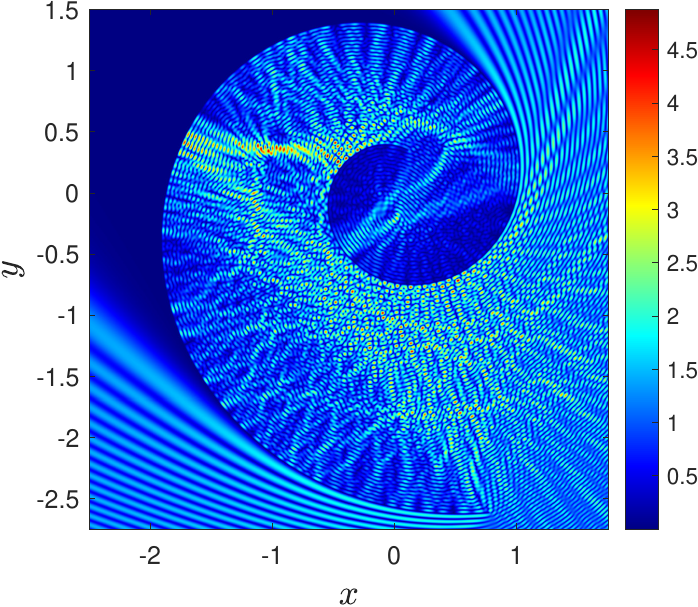}
\hspace*{2mm}
\includegraphics[height=45mm]{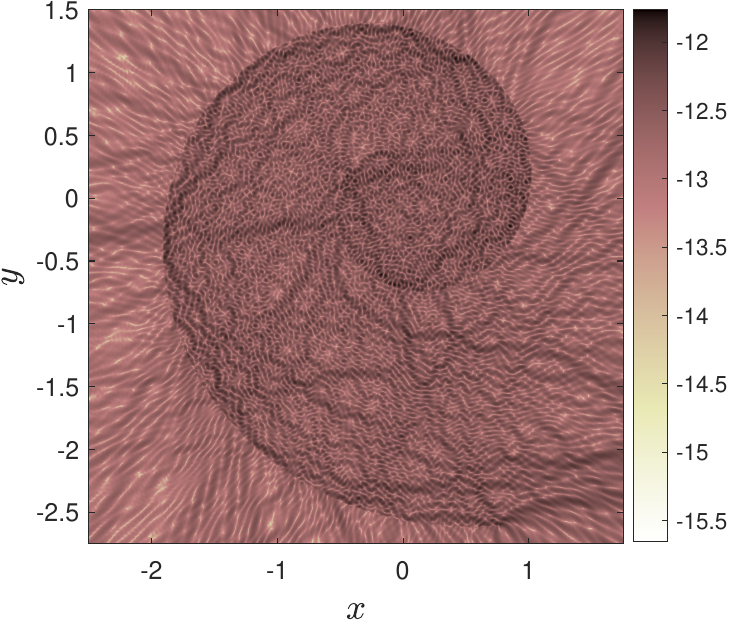}\\

\vspace{2mm}

\includegraphics[height=45mm]{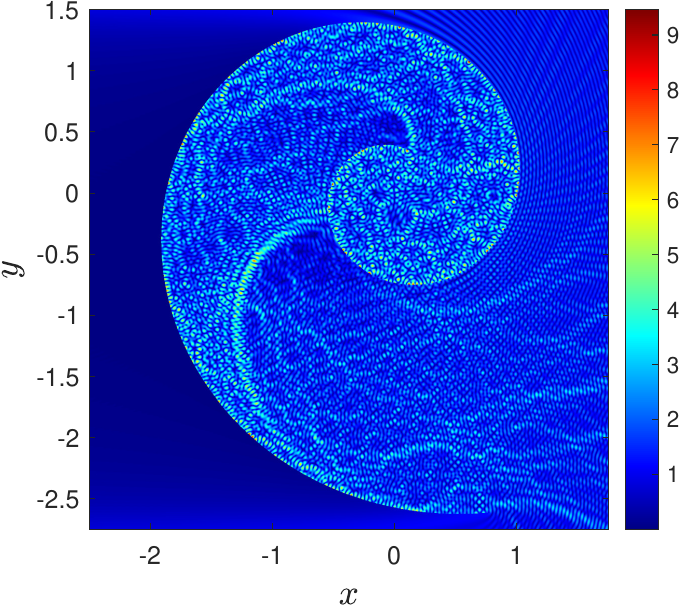}
\hspace*{2mm}
\includegraphics[height=45mm]{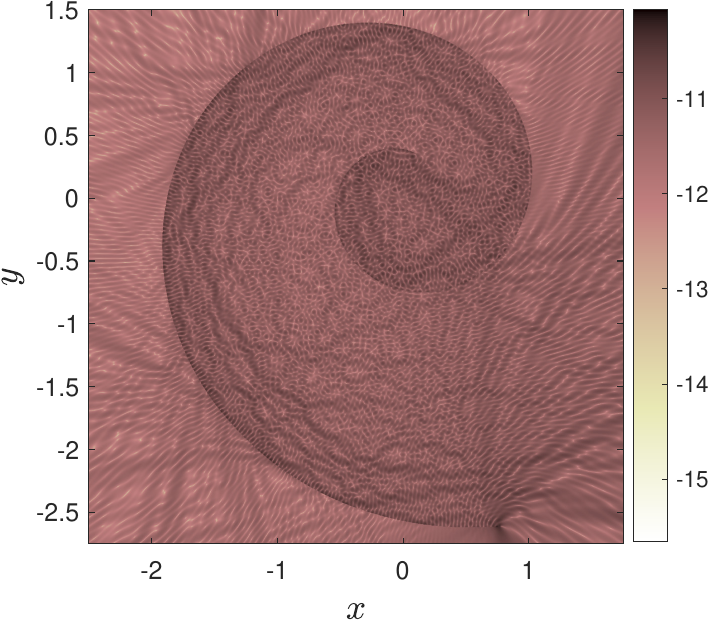}\\
\caption{\sf The magnitude of the total field $|u^{\rm tot}|$
  and the estimate abolute error for the spiral
  with $L/\lambda =200$, computed on a $3000\times 3000$ tensor grid over the rectangle
  $[-2.5, 1.75]\times [-2.75, 1.5]$. 
  Top: Dirichlet condition. Bottom: Neumann condition. Left: $|u^{\rm tot}(r)|$.
Right: $\log_{10}$ of the estimated absolute pointwise error in $u^{\rm tot}(r)$.
}
\label{fig:spiral_field}
\end{figure}

\subsubsection{Spiral}

Next, $\Gamma$ is the logarithmic spiral in \cite{BrunLint12}
\begin{equation}
r(s)=(e^s\cos(5 s),e^s\sin(5 s))\,, \quad s\in[-1,1]\,.  
\label{eq:spiral}
\end{equation}  

Figure~\ref{fig:spiral_iter} shows the conditioning and timings of
our solver for the Helmholtz Dirichlet and Neumann problems on the
spiral for $L/\lambda$ up to $25600$, compare
\cite[Table~3]{BrunLint12} for the Dirichlet problem and
\cite[Table~5]{BrunLint12} for the Neumann problem.
The relative $l_2$ error is measured on the square  $[-2.5, 1.75]\times [-2.75, 1.5]$.
The largest
relative $l_2$ errors are about $3\times 10^{-6}$ for the Dirichlet problem
and $2.7\times 10^{-6}$ for the Neumann problem.
The total number
of discretization points is roughly twice of that used
in~\cite{BrunLint12}. Similar to \cite{BrunLint12}, the number of
iterations now increases as the wavenumber increases. For
$L/\lambda=50, 200, 800$, $N_{\rm iter}=51, 60, 74$ for the
Dirichlet problem (compare $N_{\rm iter}=46, 62, 79$ in
\cite[Table~3]{BrunLint12}), and $51, 61, 74$ for the Neumann
problem (compare $48, 63, 83$ in \cite[Table~5]{BrunLint12}).
Furthermore, we observe that $T_{\rm solve}$ is very close to
$T_{\rm total}$, that is, the total computation time is dominated by
the GMRES solve time. For $L/\lambda = 800$, $T_{\rm total}\approx 17$
seconds for both the Dirichlet and Neumann problems with our
scheme, while it takes about $160$ seconds in~\cite{BrunLint12}.

Finally, Figure~\ref{fig:spiral_field}
presents high-resolution field images for the Dirichlet and Neumann problems
with $L/\lambda=200$. The incident field $u^{\rm in}$ is defined by
\eqref{eq:incfield},
using an incident angle $\theta=3\pi/4$ for the
Dirichlet problem and $\theta=\pi$ for the Neumann problem.
The boundary is divided into $440$
equi-sized panels using arc-length parameterization, resulting in $7040$
discretization points. GMRES converges to the prescribed tolerance of
$10^{-13}$ in $81$ iterations for the Dirichlet problem and $79$ iterations
for the Neumann problem.
The maximum absolute field
error on a $3000\times 3000$ grid is estimated to be less than $10^{-11}$.
Note that the image for the Dirichlet problem is directly comparable with
the top image in \cite[Figure~2]{BrunLint12}.

\subsection{Corners}
\label{sec:corners}

\begin{figure}[t]
\centering
\includegraphics[height=45mm]{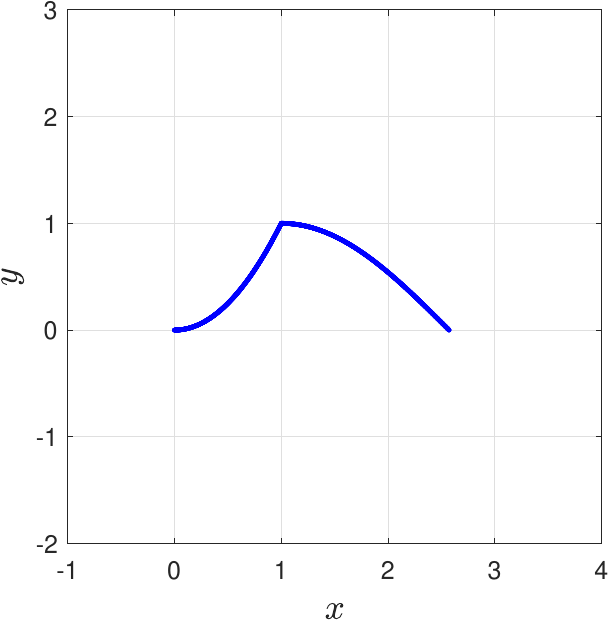}
\hspace*{2mm}
\includegraphics[height=45mm]{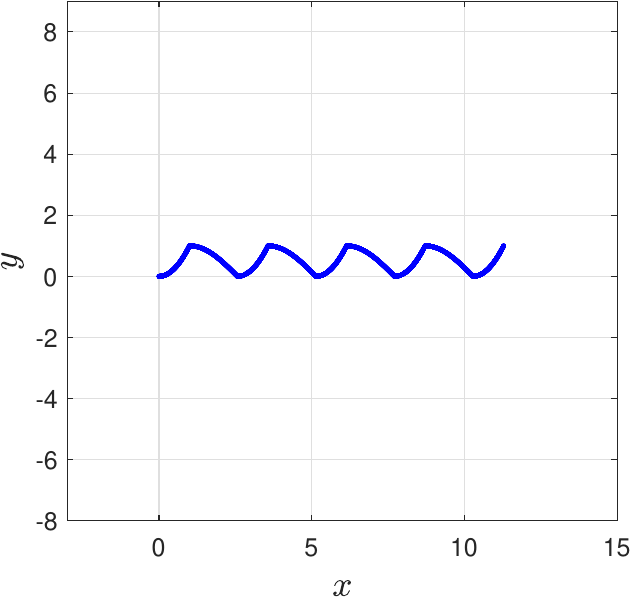}\\
\caption{\sf Boundary with corners. Left: one-corner curve
  defined by \eqref{eq:corner}. Right: eight-corner curve by tiling the curve
  on the left four times horizontally.
}
\label{fig:corners}
\end{figure}

We now consider boundaries with corners. The first boundary
consists of a parabola and a cosine curve
\begin{equation}
r(s)=
\begin{cases}
  (s, s^2)\,,\quad s\in[0,1]\,,\\
  \left(\frac{\pi s}{2},\cos(\frac{\pi s}{2})\right)\,,\quad s\in[0,1]\,,
\end{cases}
\label{eq:corner}
\end{equation}
which meet at the point $(1,1)$. We also study boundaries with
multiple corners. The left image of Figure~\ref{fig:corners} shows
the one-corner curve \eqref{eq:corner}. The right image shows an
eight-corner curve obtained by horizontally tiling this curve four
times.

\begin{figure}[t]
\centering
\includegraphics[height=45mm]{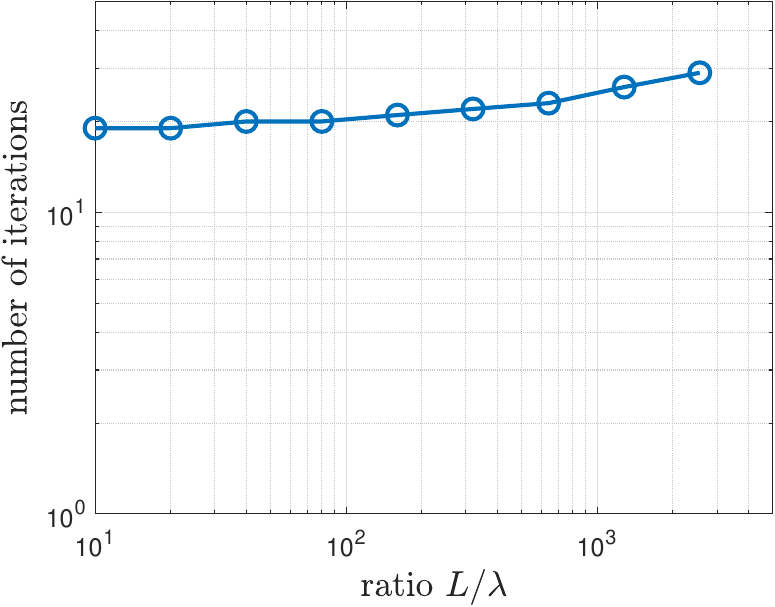}
\hspace*{2mm}
\includegraphics[height=45mm]{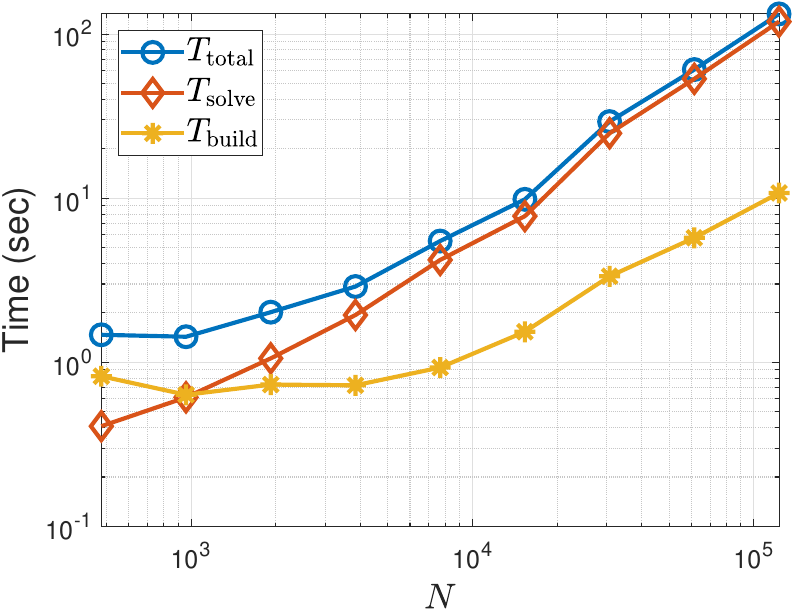}\\

\vspace{2mm}

\includegraphics[height=45mm]{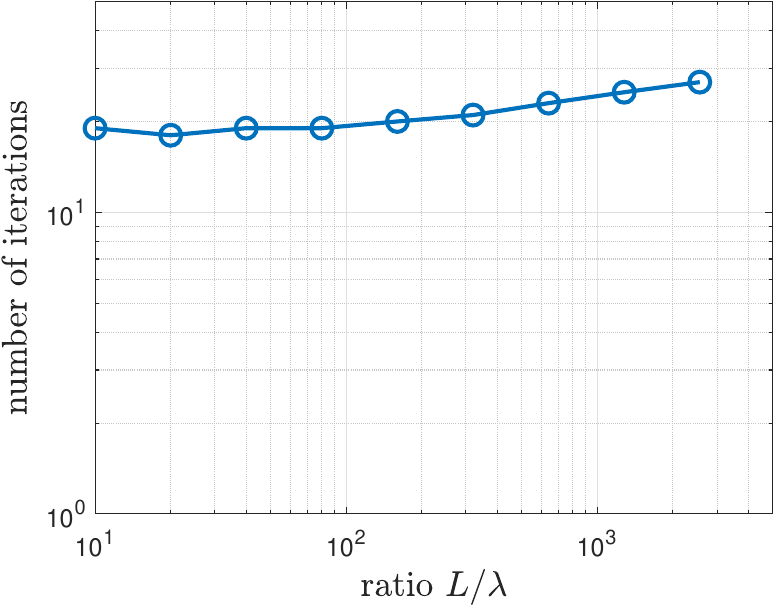}
\hspace*{2mm}
\includegraphics[height=45mm]{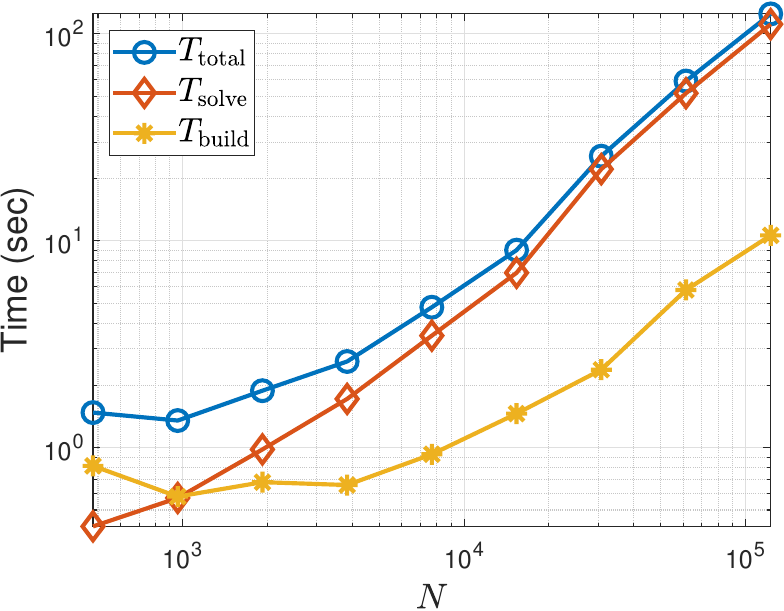}\\
\caption{\sf Same as Figure~\ref{fig:line_iter}, but for 
the one-corner curve defined by \eqref{eq:corner}.}
\label{fig:corner_iter}
\end{figure}

\begin{figure}[!ht]
\centering
\includegraphics[height=45mm]{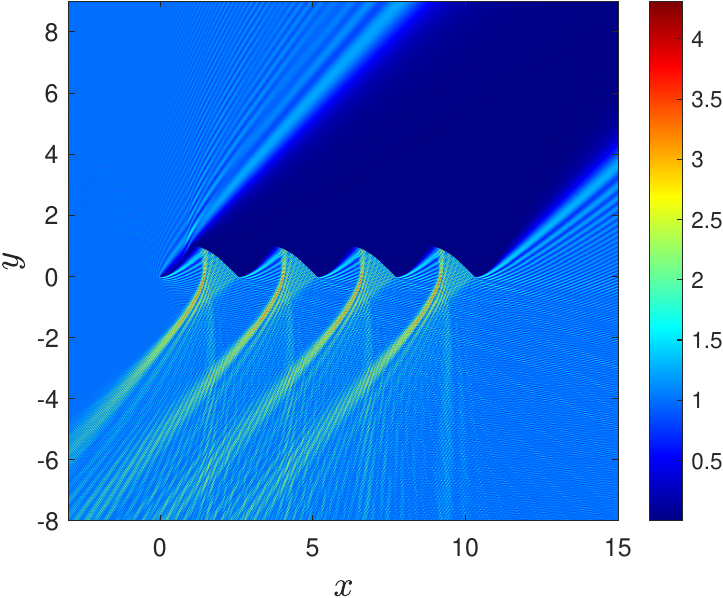}
\hspace*{2mm}
\includegraphics[height=45mm]{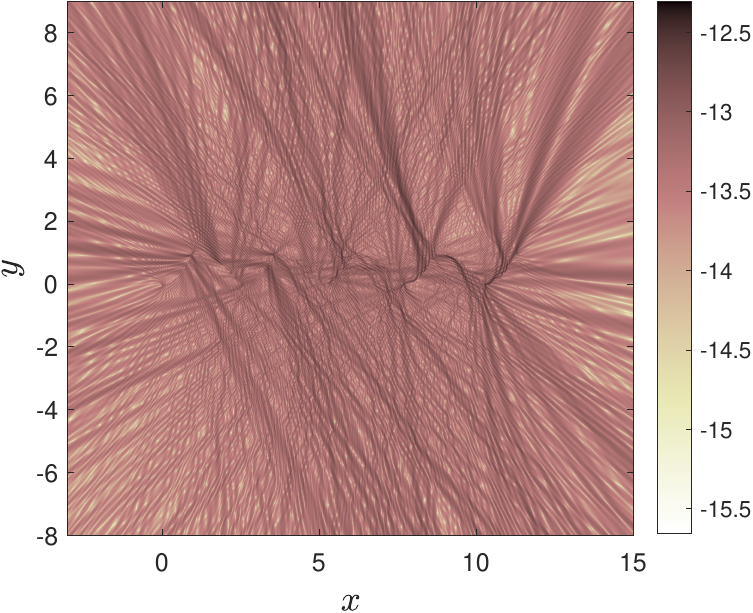}\\

\vspace{2mm}

\includegraphics[height=45mm]{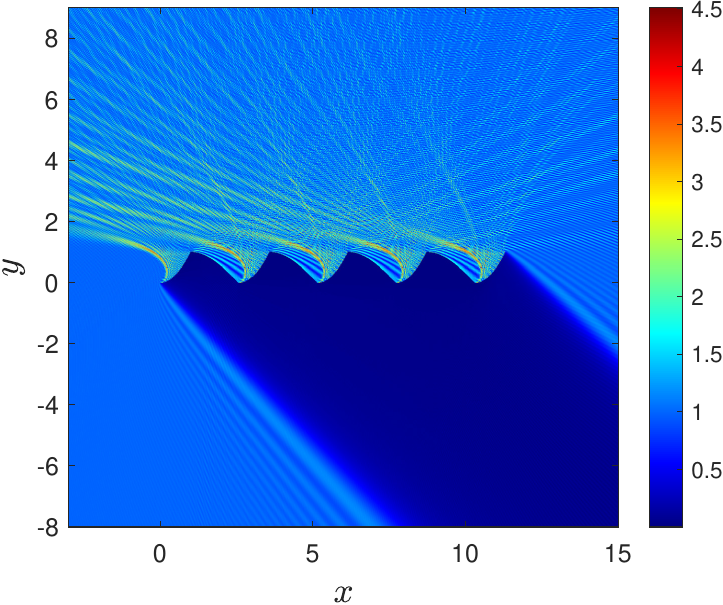}
\hspace*{2mm}
\includegraphics[height=45mm]{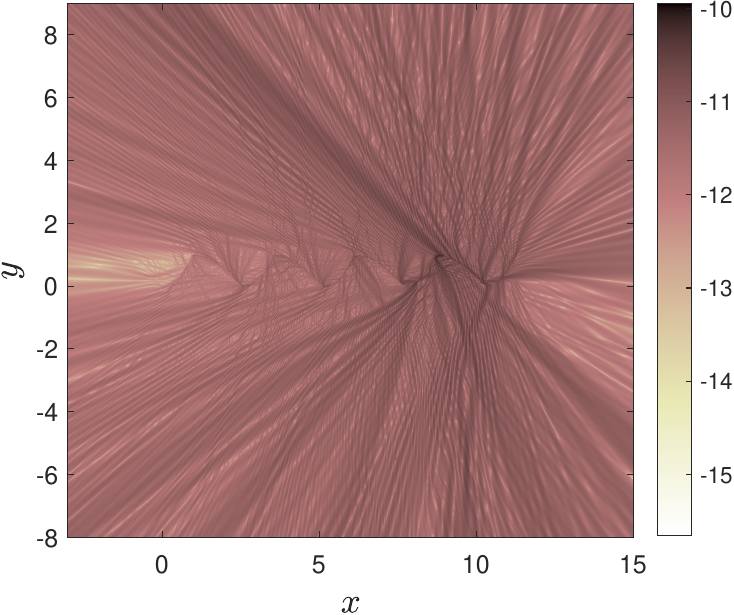}\\
\caption{\sf Same as Figure~\ref{fig:spiral_field}, but 
for the eight-corner curve.
}
\label{fig:corner_field}
\end{figure}

Figure~\ref{fig:corner_iter} presents the conditioning and timings of
our solver for the Helmholtz Dirichlet and Neumann problems on the
one-corner boundary, with $L/\lambda$ up to $2560$.
The GMRES tolerance is set to $10^{-12}$.
The relative $l_2$ error is measured on a
$300\times 300$ tensor grid over the square $[-1.5, 4]\times [-2, 3.5]$.
The largest   
relative $l_2$ errors are about $6.4\times 10^{-13}$ for the Dirichlet problem
and $1.2\times 10^{-10}$ for the Neumann problem.
We observe that the number of GMRES iterations increases mildly as the
wavenumber increases, and that the timing is again dominated by $T_{\rm solve}$.
This indicates that our numerical scheme handles corners
as effectively as endpoints.

Figure~\ref{fig:corner_field}
presents high-resolution field images for the Dirichlet and Neumann problems on the
eight-corner curve depicted in the right image of Figure~\ref{fig:corners},
with $L/\lambda=200$ (corresponding to a wavenumber $k=83.58017152213590$).
The images are computed on the rectangle $[-3, 15]\times [-8, 9]$.
The incident field $u^{\rm in}$ is defined by
\eqref{eq:incfield}, using an incident angle $\theta=\pi/4$ for the
Dirichlet problem and $\theta=-\pi/4$ for the Neumann problem.
The boundary is discretized into $9584$ points, and the GMRES algorithm converges
in approximately $40$ iterations for both problems.

\subsection{Branch points}
\label{sec:branches}

\begin{figure}[!ht]
\centering
\includegraphics[height=45mm]{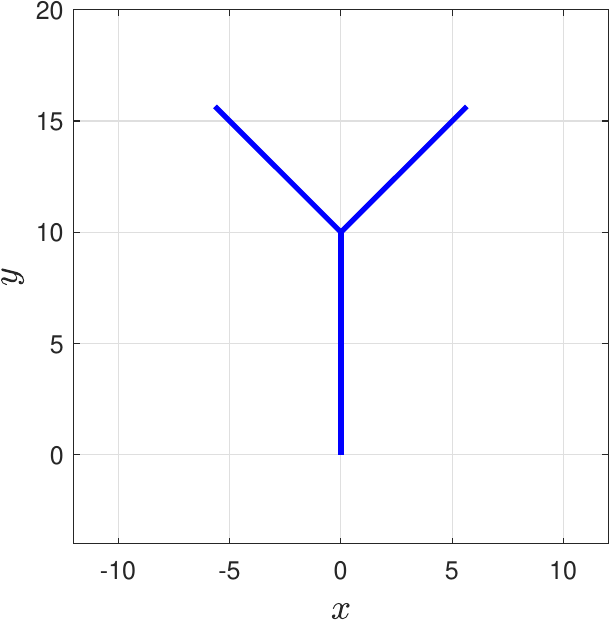}
\hspace*{2mm}
\includegraphics[height=45mm]{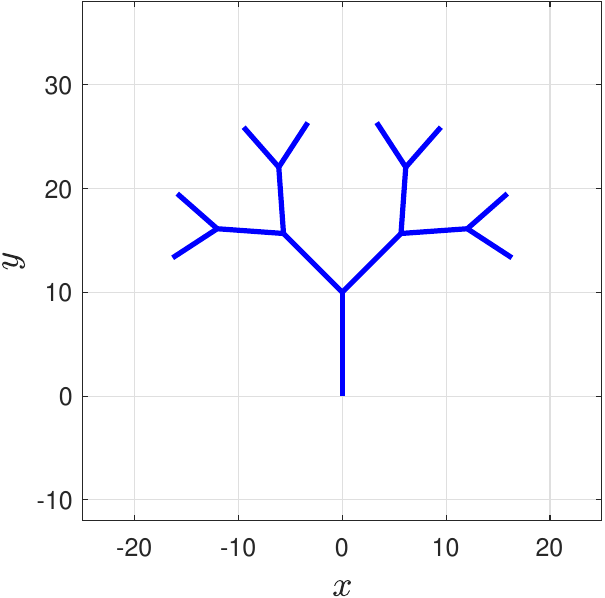}\\
\caption{\sf Branched boundary. Left: a Y shape.
  Right: a seven-branch curve with nine endpoints.
}
\label{fig:branches}
\end{figure}

We now consider boundaries with branch points. The left image of
Figure~\ref{fig:branches} shows a
Y-shaped structure composed of three line segments with lengths 
$10$, $8$, and $8$, and a top opening angle of $\pi/2$.
The ``tree'' depicted in the right image is constructed recursively by reducing
the branch lengths by a factor of $0.8$ and decreasing the opening angles by dividing by
$1.1$ at each subsequent level.

\begin{figure}[t]
\centering
\includegraphics[height=45mm]{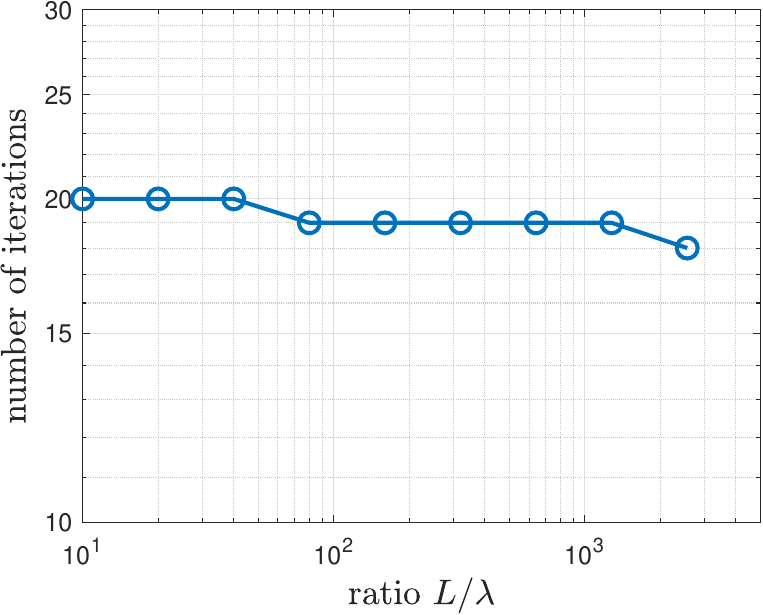}
\hspace*{2mm}
\includegraphics[height=45mm]{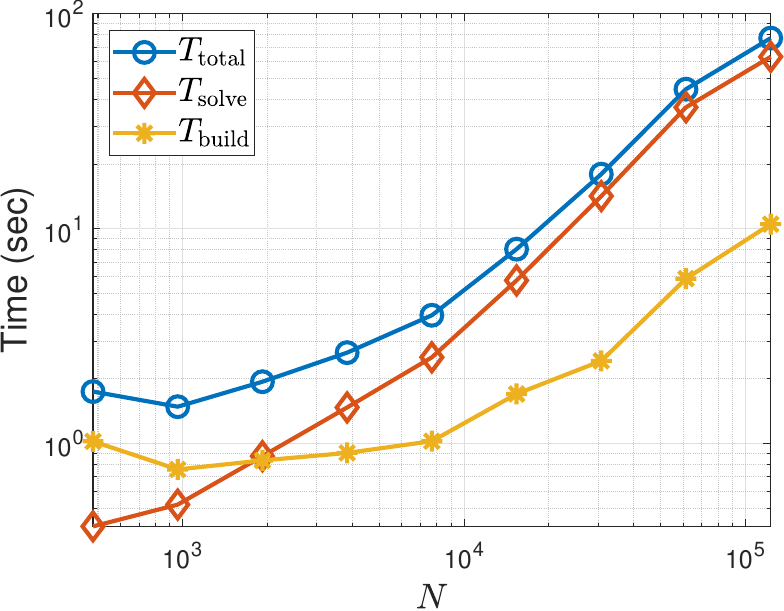}\\

\vspace{2mm}

\includegraphics[height=45mm]{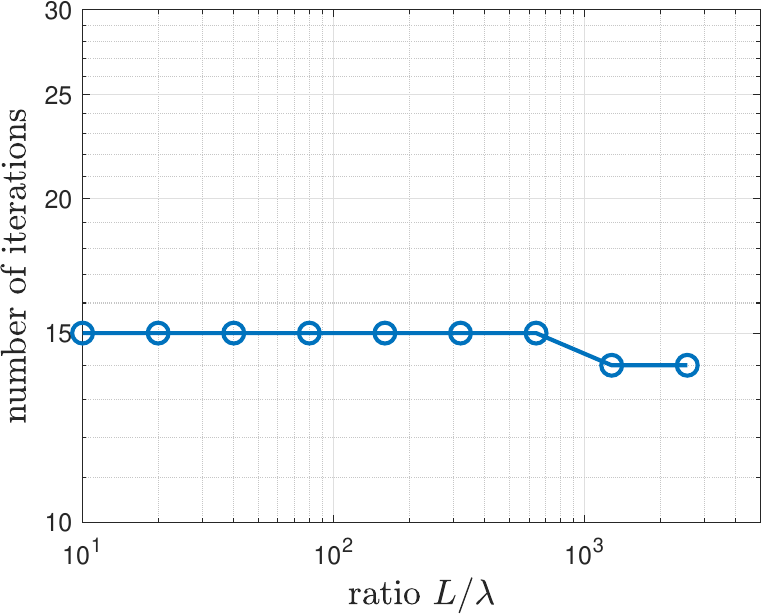}
\hspace*{2mm}
\includegraphics[height=45mm]{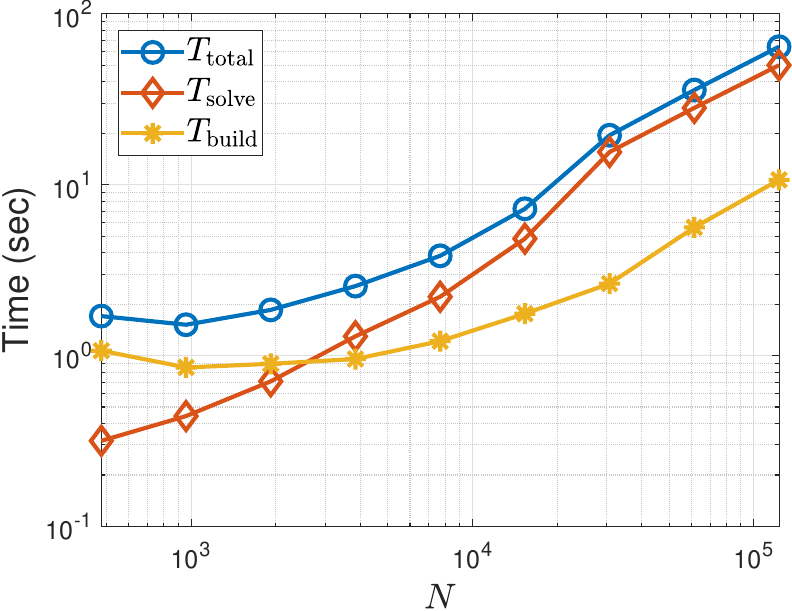}\\
\caption{\sf Same as Figure~\ref{fig:line_iter}, but
for the Y shape.}
\label{fig:branch_iter}
\end{figure}

\begin{figure}[!ht]
\centering
\includegraphics[height=45mm]{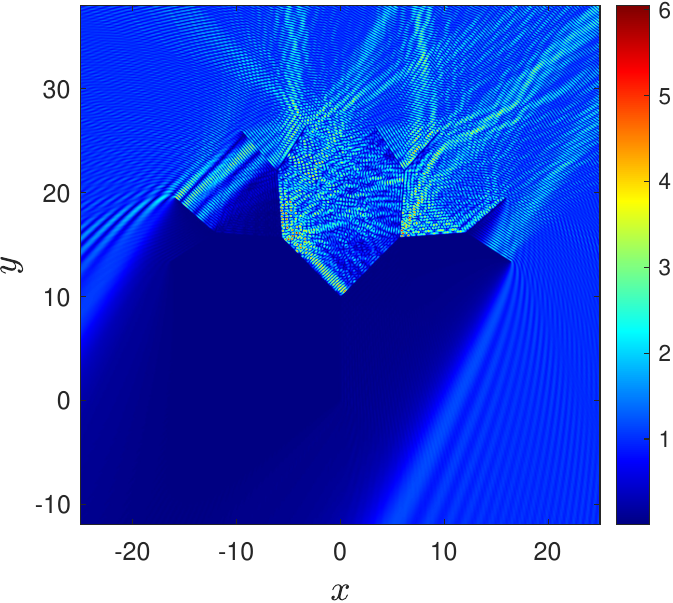}
\hspace*{2mm}
\includegraphics[height=45mm]{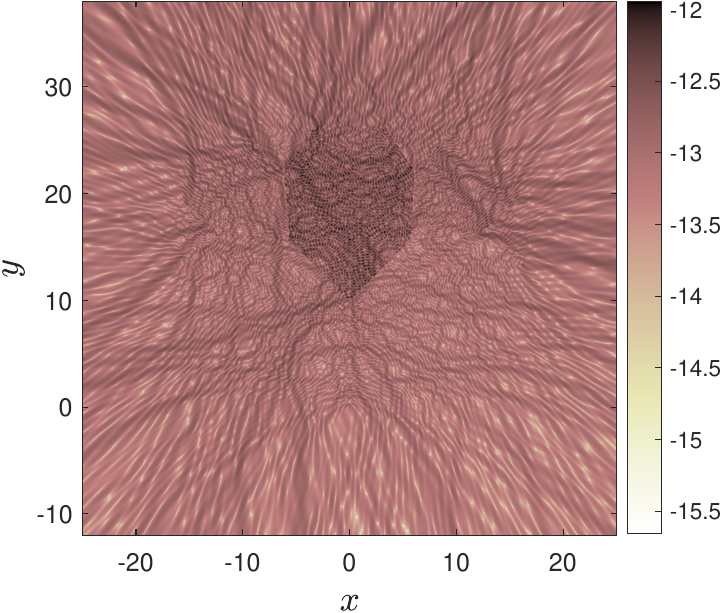}\\

\vspace{2mm}

\includegraphics[height=45mm]{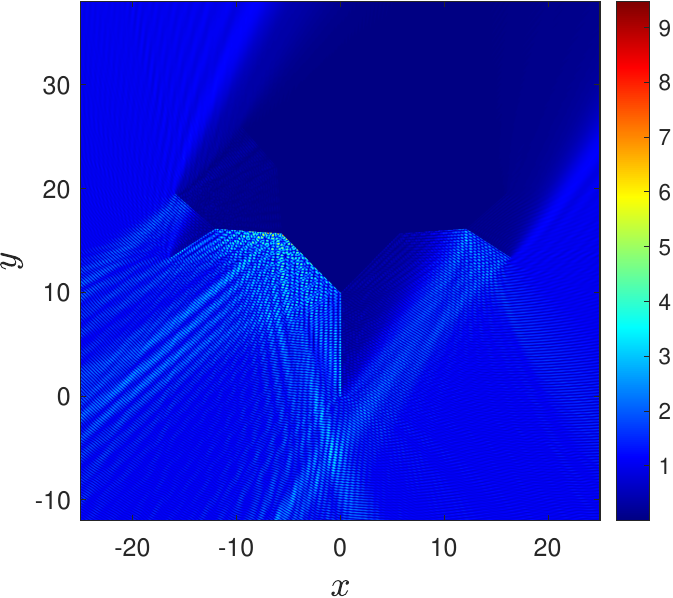}
\hspace*{2mm}
\includegraphics[height=45mm]{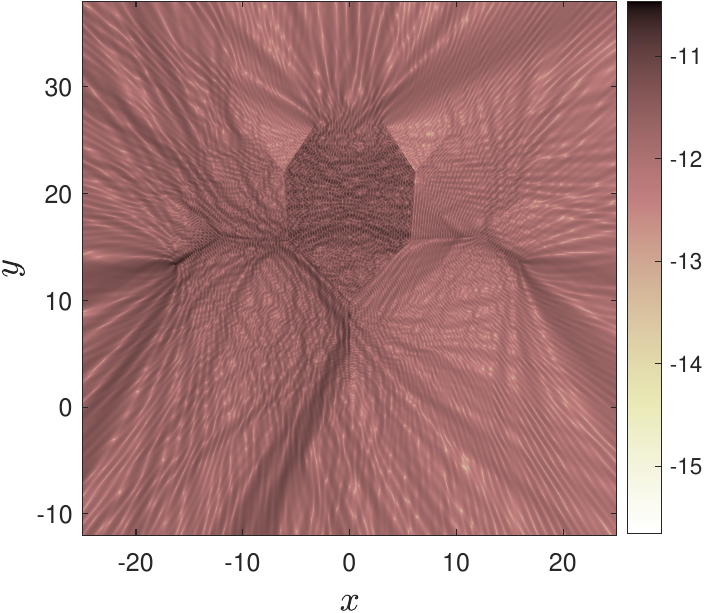}\\
\caption{\sf Same as Figure~\ref{fig:spiral_field}, but for the seven-branch
curve.
}
\label{fig:branch_field}
\end{figure}

Figure~\ref{fig:branch_iter}
illustrates the conditioning and computation times of our solver
applied to the Helmholtz Dirichlet and Neumann problems on the Y shape
depicted in the left image of Figure~\ref{fig:branches}.
The relative length
$L/\lambda$ ranges from $10$ to $2560$. The GMRES tolerance is set to $10^{-12}$.
The relative $l_2$ error is measured on a $300\times 300$ tensor grid over the square
$[-12, 12]\times [-3, 21]$. The largest relative $l_2$ errors are about 
$6.8\times 10^{-13}$ for the Dirichlet problem
and $6.4\times 10^{-11}$
for the Neumann problem.
Similar to the line segment case, we observe that the number of GMRES iterations
remains very stable as the wavenumber increases. 
Again, the computation time is dominated by $T_{\rm solve}$, demonstrating that our
numerical scheme handles branching structures as effectively as endpoints.

Figure~\ref{fig:branch_field}
presents high-resolution field images for the Dirichlet and Neumann problems
on the seven-branch curve depicted in the right image of Figure~\ref{fig:branches},
with $L/\lambda=200$ (corresponding to a wavenumber $k=13.57645917713826$). 
The images are computed on the square $[-25, 25]\times [-12, 38]$. 
The incident field $u^{\rm in}$ is defined by
\eqref{eq:incfield}, using an incident angle $\theta=-2\pi/3$ for the
Dirichlet problem and $\theta=\pi/3$ for the Neumann problem.
The boundary is discretized into $9552$ points,
and GMRES converges in $195$ iterations for the Dirichlet problem and $193$ iterations
for the Neumann problem.

\subsection{A maze}
\begin{figure}[t]
\centering
\includegraphics[height=45mm]{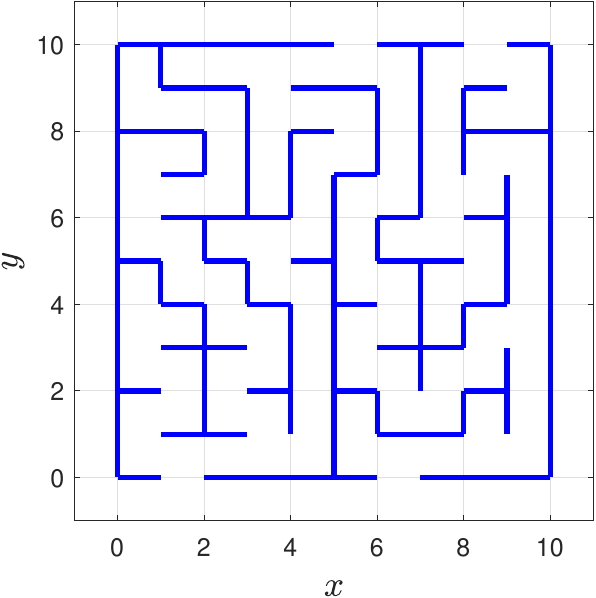}
\caption{\sf A maze-like boundary.
}
\label{fig:maze}
\end{figure}

Here, the boundary is a $10\times 10$ maze shown in Figure~\ref{fig:maze}.
It contains $30$ endpoints, $31$ corners, $18$ triple junctions,
and $2$ quadruple junctions.

\begin{figure}[!ht]
\centering
\includegraphics[height=45mm]{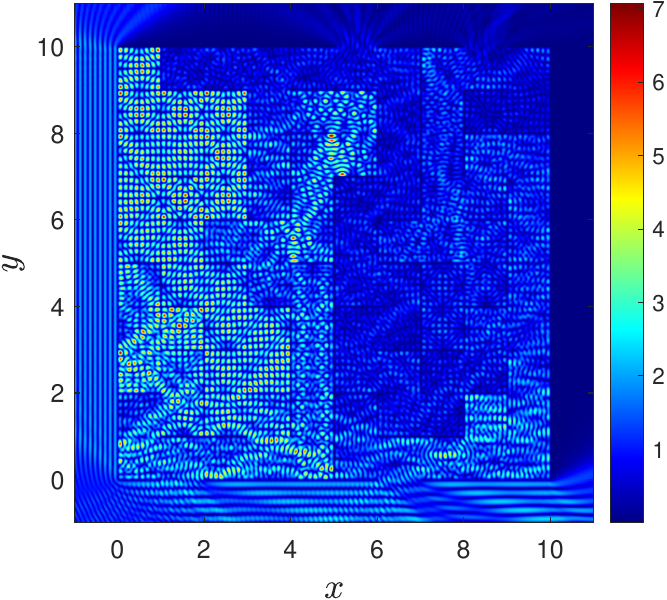}
\hspace*{2mm}
\includegraphics[height=45mm]{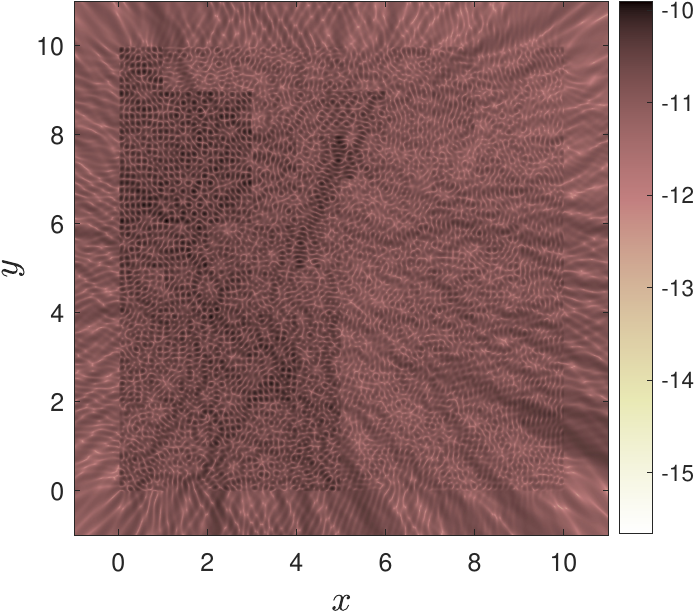}

\vspace{2mm}

\includegraphics[height=45mm]{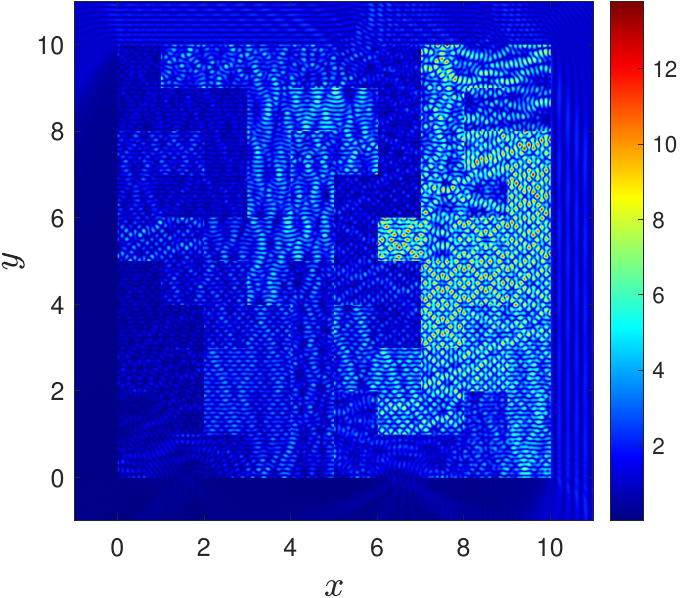}
\hspace*{2mm}
\includegraphics[height=45mm]{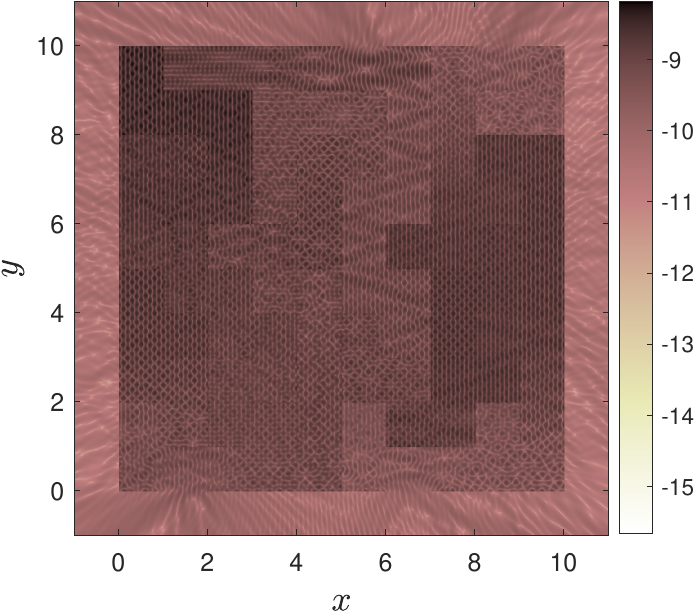}
\caption{\sf Same as Figure~\ref{fig:spiral_field}, but for the maze.
}
\label{fig:maze_field}
\end{figure}

Figure~\ref{fig:maze_field}
presents high-resolution field images for the Dirichlet problem
on the maze
with $k=10\pi$.
The incident field $u^{\rm in}$ is defined by
\eqref{eq:incfield}, using an incident angle $\theta=\pi/6$ for the Dirichlet problem
and $\theta=-2\pi/3$ for the Neumann problem. The images are computed on
the square $[-1, 11]\times [-1, 11]$.
The boundary is discretized into $28080$ points.
GMRES converges in $1877$ iterations for the Dirichlet problem
and $1815$ iterations for the Neumann problem.

\section{Conclusions and further discussions}
We have constructed an integral-equation based numerical scheme for
planar Helmholtz Dirichlet and Neumann problems on piecewise smooth
open curves. The accuracy produced by the scheme for medium-sized
and medium-wavenumber problems is, roughly speaking, $11$ digits
with Dirichlet conditions and $10$ digits with Neumann conditions.
The RCIP method is a critical component of the scheme for handling
the singularity of the unknown density in the BIE in the vicinity of
singular points of arbitrary type, including endpoints, corners, and
branch points.
Since it is used to identify the correct function spaces of the
unknown density, leading to stability and well-conditioning
comparable to that obtained in~\cite{BrunLint12} for smooth open
curves, RCIP can, in a sense, be seen as a part of our integral
equation formulation itself.

With the incorporation of the wideband fast multipole method, the
computational cost is $O(N\log N)$ for the matrix-vector product in each
GMRES iteration and $O((N+M)\log(N+M))$ for the evaluation of the field,
where $N$ is the number of source points on the boundary and $M$ is the number
of target points in the computational domain. As to the number of GMRES iterations,
it depends on the geometry and wavenumber in a highly nontrivial way 
(see, for example, \cite{marchand2022acom}). However, it is straightforward
to replace the FMM with fast direct
solvers~\cite{gopal2020acom,ho2014cpam,minden2016mms,minden2017mms,sushnikova2023mms,ye2020simaa} to deal with frequency-scan problems, scattering with multiple
right hand sides, and problems requiring large number of GMRES iterations.
Some of these solvers can be investigated using excellent software packages
such as {\tt FLAM}~\cite{ho2016flam} and {\tt chunkie}~\cite{chunkie}. Furthermore,
we observe mild loss of accuracy in our scheme for some problems
involving difficult geometries. These issues are currently under investigation.

The work in \cite{BrunLint12} has been extended to elastic scattering
in two dimensions~\cite{bruno2021ijnme} recently and to three-dimensional acoustic and
electromagnetic scattering in a preliminary report~\cite{turc2011njit}. It is
evident that our scheme can be extended to treat elastic scattering
by piecewise smooth open curves in two dimensions. However, developing high-order
efficient numerical schemes for open surface problems in three dimensions is rather
challenging. The high-order quadrature developed in \cite{zhu2022sisc} can be used to
discretize associated layer potentials in three dimensions. Nevertheless, extending the
RCIP method to handle corner and edge singularities in three dimensions remains an
open problem (see \cite{helsing2016jcp,helsing2013acha} for treatments
of axially symmetric surfaces and a cube).

\section*{Acknowledgement}
\noindent
This work was supported by the Swedish Research Council under contract
2021-03720.

\bibliographystyle{abbrv}
\bibliography{rcip}

\end{document}